\UseRawInputEncoding
\documentclass[journal]{IEEEtran}
\usepackage{amsmath,amsthm}
\usepackage{amsfonts,amssymb}
\usepackage{CJK}
\usepackage{multicol}
\usepackage{multirow}
\usepackage{diagbox}
\usepackage{graphicx}
\usepackage{subfig} %Á½¸öͼƬ²¢ÁдæÔÚ
\usepackage{lineno} %ÏÔʾÐкÅ
\allowdisplaybreaks % allow the equation display in different pages
\usepackage{epstopdf} % transfer the figure with eps formation to figure with pdf formation
\usepackage{pgf,tikz}
\usepackage{tikz}
\usepackage{graphicx}     
\usepackage{subfig} 
\usepackage{subfloat}
\usepackage{cite}
\usepackage[subfigure]{tocloft}
\usetikzlibrary{arrows,shapes,snakes,shadows,positioning,automata,patterns}
\usetikzlibrary{trees,decorations.pathmorphing,decorations.markings}
\usetikzlibrary{backgrounds}
\usepackage[%dvipdfm,  %pdflatex,pdftex %ÕâÀï¾ö¶¨ÔËÐÐÎļþµÄ·½Ê½²»Í¬
pdfstartview=FitH,
CJKbookmarks=true,
bookmarksnumbered=true,
bookmarksopen=true,
colorlinks, %×¢Ê͵ô´ËÏîÔò½»²æÒýÓÃΪ²ÊÉ«±ß¿ò(½«colorlinksºÍpdfborderͬʱעÊ͵ô)
pdfborder=001,   %×¢Ê͵ô´ËÏîÔò½»²æÒýÓÃΪ²ÊÉ«±ß¿ò
linkcolor=blue,
anchorcolor=blue,
citecolor=blue
]{hyperref}
\newtheorem{theorem}{Theorem}
\newtheorem{assumption}{Assumption}

\newtheorem{lemma}{Lemma}

\newtheorem{remark}{Remark}

\ifCLASSINFOpdf
\else
\fi

\begin{document}

\title{Solving Nonsmooth Resource Allocation Problems with Feasibility Constraints through Novel Distributed Algorithms }

\author{Xiaohong Nian, Fan Li,  Dongxin Liu
	\thanks{This work is supported in part by the National Natural Science Foundation
		of China under Grant 62173347. (Corresponding author: Fan Li)}
	\thanks{X. Nian, F. Li, and D. Liu are with Key Lab of Institute of cluster unmanned
		system, School of automation, Central South University, Changsha, 410083, China  (e-mail:funnynice@csu.edu.cn).}% <-this % stops a space
}

\maketitle

\begin{abstract}
The distributed non-smooth resource allocation problem over  multi-agent   networks  is  studied in this paper, where  each agent is subject  to  globally  coupled  network  resource   constraints   and   local  feasibility  constraints  described  in  terms  of  general  convex  sets. 
To solve such a problem, two classes of novel distributed continuous-time algorithms via differential inclusions and projection operators are proposed. Moreover, the convergence of the algorithms is analyzed by the Lyapunov functional theory and nonsmooth analysis. We illustrate that the first algorithm can globally converge  to the exact optimum of the problem when the interaction digraph is weight-balanced and the local cost functions being strongly convex. Furthermore, the fully distributed implementation of the algorithm is studied  over connected undirected graphs with strictly convex local cost functions. In addition, to improve the drawback of the first algorithm that requires initialization, we design the second algorithm which can be implemented without initialization to achieve global convergence to the optimal solution over connected  undirected graphs with  strongly convex cost functions.  Finally, several numerical simulations verify the results.
\end{abstract}


\IEEEpeerreviewmaketitle

\section{Introduction}
As the scale of the systems in practical problems such as UAV formations \cite{2016Energy}, robotic networks \cite{robost1,robost2}, sensor networks \cite{sen1} and power systems \cite{grad1,grad2,grad3} becomes increasingly huge, the traditional centralized algorithms to deal with optimization problems of large-scale network systems are not satisfactory.
Therefore, distributed algorithms that do not require a central node and can effectively reduce communication burden, are gradually receiving attention from diverse communities.

In general, according to the optimization objective, distributed optimization problems where each agent is only allowed to exchange local information with its neighbors can be classified into the two major categories. 
The optimization objective of the first category is to optimize the sum of local objective functions based on the common decision variables\cite{Zeng2015TAC,Liang2017TAC,Weiyue2020TAC,Liwei2021TAC,Kia2014auto}.
The problem considered in this paper is the other type that requires all agents to collaboratively seek the optimum for the sum of local objective functions without consensus decision variables, nevertheless the decisions of the agents are coupled with each other owing to certain coupled  constraints.  
%Due to the existence of coupled equality constraints, this class of issues is additionally known as resource allocation problems
Such problems are also known as resource allocation problems when the presence is coupled equality constraints, which have widespread applications in numerous different fields such as the economic dispach of smart grids and transportation networks\cite{2016Di,2015Dis,Xiao06,2016In,2018Dis,2018Ec}.

To cope with resource allocation problems with local feasible set constraints, continuous-time algorithms have been extensively investigated in recent years on account of their flexibility for application in real physical systems  (see \cite{Kia2016,Yi2016,Liran2020,Liang2020,Lgriad,He2017} ).
The $\epsilon$-exact penalty function is used to handle local feasible set constraints in \cite{Kia2016}. 
By combining the projection operator and the primal-dual dynamics, an initialization-free distributed algorithm was designed to solve the resource allocation problem with feasible set constraints in \cite{Yi2016}. 
Afterwards, based on \cite{Yi2016}, the work of \cite{Liran2020} simplified the way of updating auxiliary variables, thereby reducing the computational complexity.
An distributed algorithm is developed in \cite{Liang2020} with the help of singular perturbation theory and certified to converge to a suboptimal allocation of the resource allocation problem under weight-balanced digraphs.
In \cite{Lgriad}, in virtue of nonsmooth exact penalty functions to deal with local box constraints, a Laplace gradient dynamics-based algorithm is employed to address the economic dispatch problem.
Further, the work of \cite{He2017} extends the algorithm in \cite{Lgriad} to be suitable for the agents with double-integrator dynamics and illustrates through simulation examples that there is a faster convergence rate than the case with single-integrator dynamics.
It is worthwhile mentioning that the cost functions embedded in a number of engineering problems may not be differentiable (see \cite{non1,non2,non3}). Then the algorithms of the above literature  (see \cite{Kia2016,Yi2016,Liran2020,Liang2020,Lgriad,He2017}) will not be applicable, since both of them assume that the cost functions are differentiable.

Up to now, there have been a few excellent results on nonsmooth resource allocation problems (see \cite{Lian2021,Lian22021,Deng2018,Zeng2016,Zhu2019TSM}). For instance, in \cite{Lian2021,Lian22021}, distance-based exact penalty functions replace the utilization of projection operators, and an adaptive distributed algorithm is proposed for the nonsmooth resource allocation problems which have the local feasibility constraints.
In addition, for the case with nonsmooth cost functions and heterogenous local constraints, a distributed algorithm is designed in \cite{Deng2018}, via differentiated projection operators, although additional computation of the tangent cone is required and the initial state is chosen to be within the local feasible set. In \cite{Zeng2016,Zhu2019TSM}, distributed algorithms are developed by virtue of projection operators and gradient descent methods for nonsmooth resource allocation problems on undirected graphs and weight-balanced digraphs, respectively.

As we all know, it is impractical in many real applications based on an undirected communication topology among agents due to physical environment constraints and their energy limitations, and it also increases the communication costs. Additionally, there will be challenges in algorithm design and convergence analysis when encountering directed topologies, and existing algorithms such as  \cite{Kia2016,Yi2016,Liran2020,Lian2021,Lian22021,Zeng2016} may not be applicable to the resource allocation problem over directed graphs. For instance, the algorithm designed in \cite{Lgriad} can solve the case with strongly connected and weight-balanced digraphs, in which the differentiability of cost functions is indispensable.

According to the above discussions, it is evident that the nonsmooth resource allocation problem with weight-balanced digraphs and local constraints remains great research values and challenges. It is notable that the problem considered in this paper is identical to that studied in \cite{Deng2018,Deng2020,Zhu2019TSM}, but a novel distributed algorithm entirely different from others in the existing literature is developed. Moreover, a sufficient condition is given for the algorithm to be implemented in a fully distributed manner under undirected graphs. Compared to the existing literature, our main contributions have several aspects given as below.
\begin{enumerate}
	\item  We investigate the nonsmooth distributed resource allocation problem with heterogeneous feasible set constraints, which can be seem as an extension of the problems considered in \cite{2016Di,2015Dis,Xiao06,2016In,2018Dis,2018Ec,Kia2016,Yi2016,Liran2020,Liang2020,Lgriad,He2017}. Unlike \cite{Kia2016,Lgriad,He2017} where no set constraints or only box constraints are considered, the feasible set constraints here are general convex sets. Furthermore, as an improvement of \cite{Kia2016,Yi2016,Liran2020,Liang2020,Lgriad,He2017}, we consider the more general case where the cost function is nonsmooth. 
	 We develop two different classes of novel algorithms based on differential inclusions and projected output feedback for the above problem and argue their convergence resorting to  the nonsmooth analysis and the Lyapunov functional theory. Moreover, in contrast to \cite{2016Di,2015Dis,Xiao06,2016In,2018Dis,2018Ec}, the algorithms proposed in this article do not necessitate the communication of local gradient information, which is more effective in protecting privacy.
	
	%We develop two different classes of novel, in which the agent is sharing the holistic information in a sum form with its neighbors, algorithms based on differential inclusions and projected output feedback, for the above problem and analyze the convergence of
	\item We first establish the globally asymptotic convergence for the first algorithm to the exact optimal solution under a strongly connected and weight-balanced digraph with the strong convexity assumption for local cost functions. Besides, for the case of strictly convex local costs,  we characterize that the algorithm can be implemented in a fully distributed manner instead of requiring any other global information include the connectivity of the communication graph and the convexity parameters of cost functions to determine the range of the control parameters, compared to \cite{Deng2018,Zhu2019TSM,Deng2020}.
	
   \item The second algorithm can be implemented in an initialization-free way under the undirected and connected graphs without satisfying certain initial conditions of decision variables as in \cite{Yi2016,Lgriad,He2017,Lian2021,Deng2018} or auxiliary variables as in \cite{Kia2016,Lian22021,Zhu2019TSM}, which avoids the exposure of private information.
\end{enumerate}

We arrange the paper in the following order. 
Section  \ref{section2}  gives a few useful preliminaries. 
Section  \ref{section3}   formulates the nonsmooth resource allocation problem. 
Section \ref{section4}  presents the main results of this paper by designing two distributed algorithms based on projected output feedback for seeking the optimum of considered problems, and providing rigorous analysis of the convergence. 
Section  \ref{section5}  gives numerical simulations to verify the theoretical results. 
Finally, section  \ref{section6} provides final concluding remarks.

\textsl{Notations:}
$\mathbb{R}^n$ is the set of $n$-dimensional real column vectors. 
$1_n(0_n)$ denotes the $n\times1$ ones(zeros) vector.
$I_n$ is the $n \times n$ identity matrix. 
$\otimes$ stands for the kronecker product. 
$col(x_1, \ldots, x_N)=[x_1^T, \dots ,x_N^T]^T.$
$\|A\|$ and $\|x\|$ are used to represent the spectral norm of matrix $A$ and the Euclidean norm of vector $x$, respectively. 
For a set $\Omega \subset \mathbb{R}^n$, $\partial \Omega$ and $int (\Omega)$ denote the boundary points and the relative interiors of $\Omega$, respectively.
For $\Omega_1 \subset \mathbb{R}^n$ and $\Omega_2 \subset \mathbb{R}^n$, $\Omega_1 \times \Omega_2$ is utilized to denote the cartesian product.

\section{Preliminaries and Formulation}\label{section2}

\subsection{Graph Theory}
The communication graph among $N$ agents is denoted by $\mathcal{G}=\left( \mathcal{V}, \mathcal{E},\mathcal{A} \right) $, which is specified by the node set  $\mathcal{V}=\{\mathcal{V}_1,\dots ,\mathcal{V}_N\}$  the edge set $
 \mathcal{E}\in \mathcal{V}\times \mathcal{V}$ and the weighted adjacency matrix $\mathcal{A}=\left( a_{ij} \right) _{N\times N}\in \mathbb{R}^{N\times N}$.
 If the edge $ e_{ij}\in \mathcal{E} $, then $a_{ij}>0$ which indicates that node $\mathcal{V}_i$ can receive information from node $\mathcal{V}_j$; otherwise, $a_{ij}=0$.
 The path is described as a sequence of edges connecting a pair of distinct nodes. The undirected (directed) graph is connected (strongly connected) if any pair of nodes is linked by a path. The weighted in-degree and weighted out-degree of $\mathcal{V}_i$ are given by $d_{in}^{i}=\sum_{j=1}^N{a_{ij}}$  and $d_{out}^{i}=\sum_{j=1}^N{a_{ji}}$, respectively. The Laplacian matrix associated with $\mathcal{G}$  is defined as  $L=D^{in}-\mathcal{A}$, where $
 D_{in}=\text{diag}\left\{ d_{in}^{1},\dots ,d_{in}^{N} \right\} \in \mathbb{R}^{N\times N}$.
Evidently, $L1_N=0_N$.
For a connected undirected graph, $L$ is positive semidefinite and has a simple eigenvalue $0$ with the eigenvector space $\left\{ \theta \cdot 1_N| \theta \in \mathbb{R} \right\}$. Moreover, all eigenvalues of $L$ are nonnegative (see  \cite{Godsil01}). 

Besides, we can diagonalize $L$ by orthogonal transformation,  besed on the following lemma.

\begin{lemma} \label{lem.wu} (see  \cite{Yang2019TAC})
	 For a given connected undirected graph $\mathcal{G}_1$, by means of an orthogonal matrix $T=[r \; \;   R] \in \mathbb{R}^{N\times N}$, we can express $L_1$ in the following form:
\begin{align}
L_1=\left[ \begin{matrix}
	r&		R\\
\end{matrix} \right] \left[ \begin{matrix}
	0&		\\
	&		\mathcal{J}_1\\
\end{matrix} \right] \left[ \begin{array}{c}
	r^T\\
	R^T\\
\end{array} \right] 
\end{align}
where $\mathcal{J}_1$ is a diagonal matrix composed of all positive eigenvalues of $L_1$, that is, $\mathcal{J}_1=\text{diag}\left\{ \lambda_{2},\dots ,\lambda_{N} \right\} $ with $0= \lambda_1<\lambda _2\le  \dots  \le \lambda_N$ being the positive eigenvalues of $L_1$. Besides, there hold $
r=\frac{1}{\sqrt{N}}1_N,\ r^TR=0_{N}^{T},\ R^TR=I_{N-1}$, and $RR^T=I_N-\frac{1}{{N}}1_N1_{N}^{T}$.
\end{lemma}

By utilizing $Sym(L)$ to represent $\frac{L+L^T}{2}$, the following statements are the equivalent representations for a weight-balaned digraph.  
\begin{enumerate}
	\item $1_N^T L =0_N^T$;
	\item $Sym(L)$ is positive semidefinite;
	\item $\mathcal{G}$ is weight-balanced.
\end{enumerate}

Take ${\hat{\lambda} _1}, \ldots ,{\hat{\lambda} _N}$ with ${\hat{\lambda} _i} \le {\hat{\lambda} _j}$ for $i \le j$ as the  other
eigenvalues of $Sym(L)$. If $\mathcal{G}$ is a strongly connected digraph, then it follows that $0$ is a simple eigenvalue of $Sym(L)$ and the real part of all other eigenvalues is positive.

\begin{lemma}(see \cite{Kia2014auto})  \label{lemy}%\textsuperscript{\cite{Kia14}}
For a weight-balanced graph $\mathcal{G}_2$, similar to Lemma \ref{lem.wu}, with  an orthogonal matrix matrix $T=[r \; \;   R] \in \mathbb{R}^{N\times N}$, we can formulate $L_2$ as 
\begin{align}
L_2=\left[ \begin{matrix}
	r&		R\\
\end{matrix} \right] \left[ \begin{matrix}
	0&		\\
	&		\mathcal{J}_2\\
\end{matrix} \right] \left[ \begin{array}{c}
	r^T\\
	R^T\\
\end{array} \right]
\end{align}
where $
T=\left[ \begin{matrix}
	r&		R\\
\end{matrix} \right] 
$ is the same as Lemma \ref{lem.wu} and $\mathcal{J}_2=R^TL_2R$.
\end{lemma}

\subsection{Projection and Convex Analysis }
In this section, we pesent some concepts and properties about the projection operater and convex analysis (see\cite{Rockafellar70}).
The projection operater of $p$ on a closed convex set $\Omega \subset \mathbb{R}^n$ is defined as $P_{\varOmega}\left( p \right) =\arg\min _{q\in \varOmega}\lVert p-q\rVert  $ , where $p\in \mathbb{R}^n$.
\begin{lemma}\label{lemproj1}
	For a closed covex set $\Omega \subset \mathbb{R}^n$, we have the following inequilities:
	\begin{enumerate}
			\item $\lVert P_{\varOmega}\left( p \right) -P_{\varOmega}\left( q \right) \rVert \leq \lVert p-q \rVert ,\quad \forall p,\; q\in \mathbb{R}^n;$
		\item$
		\left< p-P_{\varOmega}\left( p \right) ,P_{\varOmega}\left( p \right) -q \right> \geq 0,\quad \forall ~p\in \mathbb{R}^n,~\forall ~q\in \varOmega.$
	\end{enumerate}
\end{lemma}
\begin{lemma} \label{lemproj2}
	For a closed convex set $\Omega \subset \mathbb{R}^n $, define a function on $\mathbb{R}^n \times \mathbb{R}^n$ as 
	$$
	V\left( p,q \right) =\frac{1}{2}\left( \lVert p-P_{\varOmega}\left( q \right) \rVert ^2-\lVert p-P_{\varOmega}\left( p \right) \rVert ^2 \right) .
	$$
	We can obtain that 
	\begin{enumerate}
		\item $
		V\left( p,q \right) \ge \frac{1}{2}\lVert P_{\varOmega}\left( p \right) -P_{\varOmega}\left( q \right) \rVert ^2;
		$
		\item $V\left( p,q \right) $ is continuously differentiable with respect to $p$. Further,  $
		\triangledown _pV\left( p,q \right) =P_{\varOmega}\left( p \right) -P_{\varOmega}\left( q \right).
		$	
	\end{enumerate}
\end{lemma}

%Let $\Omega \subset \mathbb{R}^n$ be an convex set, a function $f: \Omega \rightarrow \mathbb{R}$ is convex if 
%$\left( \mathcal{X}- \mathcal{Y} \right) ^T\left( \gamma _1-\gamma _2 \right) \ge 0,\ \forall \ \mathcal{X},\  \mathcal{Y}\in \ \varOmega \ $ where $\gamma _1 \in \partial f\left( \mathcal{X}\right)$, $\gamma _2 \in \partial f\left( \mathcal{Y}\right)$.
%A function  $f: \Omega \rightarrow \mathbb{R}$ is said to be strictly convex if $\left( \mathcal{X}- \mathcal{Y} \right) ^T\left( \gamma _1-\gamma _2 \right) > 0,\ \forall \ \mathcal{X},\  \mathcal{Y}\in \ \varOmega \ $  and $\mathcal{X} \ne  \mathcal{Y}$. And it is $\omega $-strongly convex ($\omega \ge 0$) if
%$\left( \mathcal{X}- \mathcal{Y} \right) ^T\left( \gamma _1-\gamma _2 \right) \ge \omega \lVert \mathcal{X}- \mathcal{Y} \rVert ^2,$ $ \forall \ \mathcal{X},\  \mathcal{Y}\in \ \varOmega .$ 
%The normal cone of $\Omega$ at $\mathcal{X}$ is defined as $
%N_{\varOmega}\left( \mathcal{X} \right) =\left\{  \mathcal{Y}:\ \left< \mathcal{Y},\mathcal{Z}-\mathcal{X} \right> \le 0,\ \forall \  \mathcal{Z}\in \ \varOmega \right\}$.

\subsection{Differential Inclusions and Nonsmooth Analysis}
In this section, we will introduce some concepts and propositions about nonsmooth analysis and differential inclusion systems. For more details, see \cite{Aubin1984Differential,Cortes08}.

%For $X,Y\subset \mathbb{R}^n$, a set-valued map $\mathcal{F}:X\rightrightarrows Y$ allocates each point $x\in X$ to a closed set $\mathcal{F}(x) \subset Y$. If for any $\varepsilon >0$, there exists a $\delta >0$ such that $\mathcal{F}\left( y \right) \subset \mathcal{F}\left( x \right) +B\left( 0,\epsilon \right)$, $\forall \,y\in B\left( x,\delta \right)$, then the set-valued map $\mathcal{F}$ is upper semicontinuous.

%A function $f:\mathbb{R}^n\rightarrow \mathbb{R}$ is said to be locally Lipschitz near $x\in \mathbb{R}^n$ if there exist $\theta _x,\delta >0$, such that $\|f\left( z \right) -f\left( y \right) \|\leq \theta _x\|z-y\|,\forall \,y,z\in B\left( x,\delta \right)$. 

For a locally Lipschitz function $f:\mathbb{R}^n\rightarrow \mathbb{R}$, the Clarke's generalized gradient $\partial f$ is defined as $$\partial f=co\left\{ \underset{k\rightarrow +\infty}{\lim}\nabla f\left( x_k \right) \left| x=\underset{k\rightarrow +\infty}{\lim}x_k,\ x_k\notin  \mathcal{O}\cup \varOmega _f \right. \right\}, $$ where $co\left\{ \cdot \right\} $ represents convex hull, $\varOmega _f$
denotes the set of points in which $f$ is not diffrentiable, and $\mathcal{O}\subset \mathbb{R}^n$ is a set of Lebesgue measure zero. And if $f$ is convex, then the Clarke's generalized gradient is consistent with the sub-differential. It is known that $\partial f$ takes nonempty, compact and convex values and is locally bounded and upper semicontinuous.

A differential inclusion is given by
\begin{align}\label{1}
\dot x \in \mathcal{F}(x), \quad x(0) = x_0
\end{align}
where $\mathcal{F}:\mathbb{R}^n\rightrightarrows \mathbb{R}^n$ represents a set-valued map. An absolutely continuous map $x: [0, ~ T] \rightarrow \mathbb{R}^n$ is called a Caratheodory solution of (\ref{1}) on $[0, ~ T]$, if $x$ satisfies (\ref{1}) for almost all $t\in \left[ 0,T \right] $.

%\begin{lemma}\label{map}
%If the set-valued map $\mathcal{F}$ is upper semicontinuous, locally bounded and takes nonempty compact convex values, then there is a Caratheodory solution to (\ref{1}) for any initial state.
%\end{lemma}

\begin{lemma}\label{map}
	If $\mathcal{F}$ is an upper semicontinuous and locally bounded set-value map, and it takes nonempty, compact, and convex values, then it can be concluded that there is a Caratheodory solution to (\ref{1}) for any initial state.
\end{lemma}

The set-valued Lie derivative of a continuous differentiable function $V:\mathbb{R}^n\rightarrow \mathbb{R}$ along with (\ref{1}) is defined as 
%$$
%\mathcal{L}_\mathcal{F}V = \{a \in \mathbb{R} \, | \, \exists \, v \in \mathcal{F}(x) {\text{ s.t. }} v^T \nabla V(x) =a \}.
%$$
$$
\mathcal{L}_\mathcal{F}V =\{ v^T\nabla V|v \in  \mathcal{F}(x)\}.
$$
The set-valued LaSalle invariance principle as below is essential to the subsequent proof of convergence.

\begin{lemma}\label{laser}
 Let $V:\mathbb{R}^n\rightarrow \mathbb{R}$ be a continuously differentiable function and $S\in \mathbb{R}^n$ be a compact and strongly positively invarint set for (\ref{1}). 
Assume that the Lie derivative satisfies $\max \mathcal{L}_\mathcal{F}V \leq 0$ or $\mathcal{L}_\mathcal{F}V = \emptyset$ for all  $ x \in S$, and the Caratheodory solutions of (\ref{1}) are bounded, then the solutions of (\ref{1}) with any initial point in $S$ converges to the largest weakly positively invariant set $\mathcal{M}\subset S\cap \{x\in \mathbb{R}^n\,|\,0\in \mathcal{L}_{\mathcal{F}}V\left( x \right) \}$.
\end{lemma}

\section{Problem Formulation}\label{section3}
The constrained resource allocation problem concerned in this article can be described as follows:  
\begin{align} \label{op1}
\min_{y\in \mathbb{R}^{Nn}}&~~f\left( y \right) ,\quad f\left( y \right) =\sum_{i=1}^N{f_i}\left( y_i \right) \notag \\ &\text {subject to}~ \sum_{i=1}^N y_{i} = \sum_{i=1}^N d_{i} \notag  \\&  y_{i} \in\Omega_{i}, ~~ i \in \{1,\dots,N\}
\end{align}
where $y=col(y_1,\dots,y_N) \in \mathbb{R}^{Nn} $ consists of the local decision $y_i$ satisfying the local feasible set, i.e., $y_i \in \Omega_{i}$,  and $f_i: \Omega_i \rightarrow \mathbb{R}$ is the nonsmooth local cost function. Moveover,   $\sum_{i=1}^Ny_{i} = \sum_{i=1}^Nd_{i}$  is the network  resource constraint, where $d_i \in \mathbb{R}^n$ is the local resource.

We aim to design effective distributed algorithms for the constrained problem \eqref{op1} such that each agent minimizes the global cost function while sharing private information only with  its neighbors.
\begin{remark}
The resource allocation problem considered in this paper allows the local cost function to be non-smooth, which extends the problem considered in \cite{Kia2016,Yi2016,Liran2020,Liang2020,Lgriad,He2017}. Meanwhile, the locally feasible set is a general convex set, while only the special case of the box constraint is considered in \cite{2016Di,2015Dis,Xiao06,2016In,2018Dis,2018Ec}.
\end{remark}

The following mild assumptions utilized in the subsequent analysis is meaningful and widely used in the literature (see    \cite{Yi2016,Deng2018,Deng2020,Zhu2019TSM}).

\begin{assumption}\label{as1}
(Slater's condition) For each $i \in \{1,\dots,N\}$, there exists a solution $x_i\in int\left( \varOmega _i \right) $  such that  $\sum_{i=1}^{N}x_i=\sum_{i=1}^{N}d_i$. 
\end{assumption}

\begin{assumption}\label{as2}
For each $i\in \{ 1, \dots, N\}$, $f_i$ is convex and locally Lipschitz continuous.
\end{assumption}

The following lemma is the optimality condition for problem \eqref{op1}.
\begin{lemma}\label{kkt}(see \cite[Theorem 3.34]{Franz2}) For each $i \in \{1,\dots,N\}$, 
     $y_i^*\in \Omega_{i}$ is the optimal solution of \eqref{op1}, if and only if there exists $s \in \mathbb{R}^n $ such that 
    \begin{subequations} \label{K}
    	\begin{align}   %%%%  algo1???
   0_{Nn}\in \partial f_i\left( y_{i}^{\ast} \right) -s^{\ast}+N_{\varOmega _i}\left( y_{i}^{\ast} \right)\\
\sum_{i=1}^N{y_{i}^{\ast}}=\sum_{i=1}^N{d_i}.
 \end{align} 
 \end{subequations}
\end{lemma}

\begin{lemma}\label{lem.dac}%\textsuperscript{\cite{Kia14}}
	For any matrix $A\in \mathbb{R}^{n\times n}$, and vector $x\in \mathbb{R}^n$, we have that $Ax=0_n$ if and only if $A^TAx=0_n$.
\end{lemma}
\emph{Proof:}  If there exists $\bar{x} \in \mathbb{R}^n$ such that $A^TA\bar{x}=0_n$, then we can imply that $\left( A\bar{x} \right) ^T\left( A\bar{x} \right)=\bar{x}^TA^TA\bar{x}=0$, which means that $A\bar{x}=0_n$.

Conversly, it is not hard to get $A^TA\tilde{x}=0_n$ from $A\tilde{x}=0_n$ for $\tilde{x} \in \mathbb{R}^n$.\hfill $\Box$

%The aim of these agents is to minimize the globalcost function $f\left( y \right) $?where the global cost function is the sum of local cost functions of agents, that is, $f\left( y \right) =\sum_{i=1}^N{f_i\left( y_i \right)}$ with $y=col\left( y_1,\cdots ,y_N \right) $ Moreover, the minimum satisfies thelocal constraints of all agents and the network resource constraint. Specifically, the agents face the following optimization problem:

\section{Main Result}\label{section4}
In this section, two classes of algorithms are proposed to tackle the distributed resource allocation problem (\ref{op1}) in Section \ref{A}. Afterwards, the convergence properties of  two algorithms are discussed in Section \ref{B} and Section \ref{C}, respectively.

\subsection{Distributed  Algorithm Design} \label{A}

In this section, we focus on the design of the distributed algorithms on the basis of projected output feedback for the problem (\ref{op1}).

In order to drive all agents to minimiz the global cost function while utilizing only local information, we design the distributed algorithm as follows:
\begin{align} \label{al1}
	\begin{cases}      %%%algo1
		\dot{x}_i\in y_i-x_i-\partial f_i\left( y_i \right) +s _i\\
		\dot{s}_i=k_1\left( w _i-y_i+d_i \right) +k_2\sum_{j=1}^N{a_{ij}\left( s _j-s _i \right)}\\
		\dot{w}_i=k_3\sum_{j=1}^N{a_{ij}\left( \left( y_j-d_j+w _j \right) -\left( y_i-d_i+w _i \right) \right)}\\
		y_i=P_{\varOmega _i}\left( x_i \right) 
	\end{cases}
\end{align}
where
\begin{align} \label{can1}
k_1>\frac{\lVert L \rVert ^2}{\hat{\lambda}_2\omega}, \quad k_2>\frac{k_{1}^{2}}{\hat{s}_{2}^{2}},  \quad  k_3>0.
\end{align}

\begin{remark}
In this system, $\partial f_i\left( y_i \right) $ is utilized to seek the optimum of the nonsmooth problem \eqref{op1}, $s_i$, and $w_{i}$ are the auxiliary variables and the projected output feedback term $P_{\varOmega _i}\left( x_i \right) $ is introduced to solve the set constraints. Moreover, the projected output feedback term enable the initial state $x_i(0)$ to be outside of the set constraints, which is not permitted in\cite{Yi2016,Deng2018,Zeng2016}.
\end{remark}

For the sake of relaxing  initial value demands of the algorithm (\ref{al1}) for auxiliary variables, we next introduce an initialization-free algorithm as follows:
\begin{align} \label{al2}
	\begin{cases}        %%%algo2
		\dot{x}_i\in y_i-x_i-\partial f_i\left( y_i \right) +s _i\\
		\dot{s}_i=k_1\left( w _{i}^{L}-y_i+d_i \right) +k_2\sum_{j=1}^N{a_{ij}\left( s _j-s _i \right)}\\
		\dot{w}_i=k_3\sum_{j=1}^N{a_{ij}\left( \left( y_j-d_j+w _{j}^{L} \right) -\left( y_i-d_i+w _{i}^{L} \right) \right)}\\
		y_i=P_{\varOmega _i}\left( x_i \right) 
	\end{cases}
\end{align}
where\ $w _{i}^{L}=\sum_{j=1}^N{a_{ij}\left( w _i-w _j \right)}$,
\begin{align} \label{can2}
k_1>\frac{\lVert L \rVert ^2}{{\lambda}_2^2\omega}, \quad
k_2>\frac{k_{1}^{2}\lVert L \rVert ^2}{{\lambda}_2^3} ,  \quad  k_3>0.
\end{align}
\begin{remark}
	The structure of algorithm (\ref{al2}) is similar to  algorithm (\ref{al1}), the main difference being that the auxiliary variable $w _{i}^{L}$ is replaced by the $w _{i}$, which allows the algorithm to be implemented in an initialization-free manner.
\end{remark}
\begin{remark}
Notice that according to the aforesaid two algorithms (\ref{al1}) and (\ref{al2}), what information really transmitted between the agent and its neighbors is actually the overall information $y_j-d_j+w_j(w_j^L)$ in the form of a sum, rather than the specific private information $y_j$, $d_j$, or $w_j(w_j^L)$ alone. In this case, the agent has no way to identify specific private information from the overall information $y_j-d_j+w_j(w_j^L)$, and the neighbor's $y_j$, $d_j$, and $w_j(w_j^L)$ are still unknown to the intelligence, which avoids privacy leakage in a certain sense.
Moreoer, based on such a way to share sum information, it does not increase the additional communication burden of the network, campared to the current algorithms in \cite{Yi2016,Deng2018,Zhu2019TSM,Deng2020}.
\end{remark}
\begin{remark}
Through the analysis in the sequel, similar to \cite{Deng2018,Deng2020,Zhu2021}, we will see that the convergence of the algorithm requires the parameters $k_1$ and $k_2$ to meet a specific range, which depends on $
\hat{\lambda}_2, \, \lVert L \rVert \ $ and $\omega $.  We can pre-calculate the dependent values through an additional distributed consensus algorithm (see \cite{concensus}) to determine the range of the parameters $k_1$ and $k_2$  in advance.
\end{remark}

\subsection{Convergence Analysis Of Algorithm \eqref{al1}} \label{B}
In this section, we analyze the characteristics of the equilibrims and the convergence of \eqref{al1}. Specifically, the proof of the convergence is established on the nonsmooth analysis and the Lyapunov functional theory.

Let 
$x=col( x_1,x_2,\ldots ,x_N ),\; 
s =col( s _1,s _2,\ldots ,s _N ), \;	
w =col( w _1, w _2,\ldots ,w _N ),\;
d=col( d_1,d_2,\ldots ,d_N ), \;
\varOmega =\varOmega _1 \times \varOmega _2 \times \ldots \times \varOmega _N , \;
\partial f( y ) =col( \partial f_1( y_1 ) ,\partial f_2( y_2 ) , \ldots ,\partial f_N( y_N ) ), \;
y=col( y_1,y_2,\ldots ,y_N ) 
$.

 We can recast (\ref{al1}) in a campact form as:
\begin{align} \label{al11}
	\begin{cases}          %%%%algo1  ????
		\dot{x}\in y-x-\partial f\left( y \right) +s \\
		\dot{s}=k_1\left( w -y+d \right) -k_2\left( L\otimes I_n \right) s \\
		\dot{w}=-k_3\left( L\otimes I_n \right) \left( w -y+d \right) \\
		y=P_{\varOmega}\left( x \right) .
	\end{cases}
\end{align}

Natably, the existence of the solution of (\ref{al11}) can be guaranteed by Lemma \ref{map} under Assumption \ref{as2}.

We have the following result with regard to the equilibrium of (\ref{al11}).
\begin{theorem}\label{the1}
For the nonsmooth resource allocation problem \eqref{op1}, consider the case where the communication topology is a strongly connected and weight-balanced digraphs.
 If Assumptions \ref{as1} and \ref{as2} hold, with the initial condition satisfying $\sum_{i=1}^N{w_i(0)}=0_n$,  then $(y^*,x^*,s^*,w^*)$ is an equilibrium point of \eqref{al11}, if and only if $y^*$ is an optimal solution of \eqref{op1}.
\end{theorem}

\emph{Proof:}
1) We assume that $(y^*,x^*,s^*,w^*)$ is an equilibrium  of (\ref{al11}), then one has 
\begin{subequations} \label{m}
	\begin{align}\label{Th1a}    %%%%  algo1???
		0_{Nn}\in& y^{\ast}-x^{\ast}-\partial f\left( y^{\ast} \right) +s ^{\ast} \\  \label{Th1b}
		0_{Nn}=&k_1\left( w ^{\ast}-y^{\ast}+d \right) -k_2\left( L\otimes I_n \right) s ^{\ast}\\ \label{Th1c}
		0_{Nn}=&-k_3\left( L\otimes I_n \right) \left( w ^{\ast}+y^{\ast}-d \right) \\   \label{Th1d}
		y^{\ast}=&P_{\varOmega}\left( x^{\ast} \right).   
	\end{align} 
\end{subequations}

Note that $L1_N=0_N$ and  $1_{N}^{T}L=0_{N}^{T}$ are satisfied due to the strony connectedness and weight-balance of digraphs. Then there exists $\theta _1\in \mathbb{R}^n$ such that $\left( L\otimes I_n \right) s ^*=1_N\otimes \theta _1$ form (\ref{Th1b}) and (\ref{Th1c}). Subsequently, one can obtain that $\left( L^TL\otimes I_n \right) s ^*=L^T1_N\otimes \theta _1=0_{Nn}$,  which indicates that $\left( L\otimes I_n \right) s ^*=0_{Nn}$ i.e. $
s _{i}^{*}=s _{j}^{*},\ \forall \ i,j\in \left\{ 1,\dots,N \right\} $ by Lemma \ref{lem.dac}. Further, it results from (\ref{Th1b}) that $k_1\left( 1_{N}^{T}\otimes I_n \right) \left( w ^{\ast}-y^{\ast}+d \right) =k_2\left( 1_{N}^{T}L\otimes I_n \right) s ^{\ast}=0_{n}$. Additionally, it follows from (\ref{al11}) that $(1_N^T\otimes I_n)\dot{w}(t)=0_n,$ and applying $\sum_{i=1}^N{w _i\left( 0 \right) =0_n}$, i.e. $( 1_{N}^{T}\otimes I_n ) w ( 0 ) =0_{n}$, we know that $
(1_{N}^{T}\otimes I_n) {w}\left( t \right) =0_{n} $ for any $t \ge 0$. Therefore, one can get that $
k_1\left( 1_{N}^{T}\otimes I_n \right) \left( -y^*+d \right) =0_{n}$, i.e. $\sum_{i=1}^N{y_{i}^{*}}=\sum_{i=1}^N{d_i}$.

Besides, from (\ref{Th1a}) and (\ref{Th1d}), it can be seen that $y^{\ast}=P_{\varOmega}\left( y^{\ast}-\partial f\left( y^{\ast} \right) +s ^{\ast} \right) $, that is, $0_n\in \partial f\left( y_{i}^{*} \right) -s _{i}^{*}+N_{\varOmega _i}\left( y_{i}^{*} \right) $.
According to the above analysis,  it is concluded that $y^*$ is an optimal solution of the problem (\ref{op1}) with reference to Lemma \ref{kkt}. 

2) If $y^*$ is an optimal solution of the problem (\ref{op1}), then there exists $s_{i}^* \in \mathbb{R}^n$, such that
\begin{align}  \label{zui}
	\begin{cases}        
0_{n}\in \partial f_i\left( y_{i}^{\ast} \right) -s _{i}^{\ast}+N_{\varOmega _i}\left( y_{i}^{\ast} \right) \\
\sum_{i=1}^N{y_{i}^{\ast}}=\sum_{i=1}^N{d_i}\\
s _{i}^{\ast}=s _{j}^{\ast}\ \forall i,j\in \{1,\ldots,N\}\\
y_{i}^{\ast}\in \varOmega _i\ \forall i\in \{1,\ldots,N\}.
	\end{cases}
\end{align} 

By taking $s ^*=col\left\{ s _{1}^{*},...,s _{N}^{*} \right\}$, it is obvious that $\left( L\otimes I_n \right) s ^*=0_{Nn}$ based on (\ref{zui}).

Next, we can take $x^*=col\left\{ x_{1}^{*},...,x_{N}^{*} \right\}$, where $x_{i}^{*}\in y_{i}^{*}-\partial f_i\left( y_{i}^{*} \right) +s _{i}^{*}$. Then, from (\ref{zui}), it follows that (\ref{Th1a}) and (\ref{Th1d}) hold.

Furthermore, denote $w ^*=col\left\{ w _{1}^{*},...,w _{N}^{*} \right\}$, where $w _{i}^{*}=y_{i}^{*}-d_i$, then (\ref{Th1b}) and (\ref{Th1c}) are satisfied, since $\left( L\otimes I_n \right) s ^*=0_{Nn}$. Therefore, $\left( y^*,x^*,s ^*,w ^* \right)$ is an equilibrium  of (\ref{al11}).\hfill $\Box$

In the following, the asymptotic convergence associated with  algorithm (\ref{al1}) is studied over a weight-balanced digraph and a undirected graph.
\begin{theorem}\label{the2}
For the problem (\ref{op1}) with Assumptions \ref{as1} and \ref{as2}, consider the case where  the local cost functions are $\omega$-strongly convex and the communication topology is a strongly connected and weight-balanced digraph.
 Suppose that the  initial point $(y(0),x(0),s(0),w(0))$ satisfies $\sum_{i=1}^N{w _i\left( 0 \right) =0_n}
$, then the algorithm (\ref{al1}) can converge asymptotically to the optimum of problem (\ref{op1}).
\end{theorem}
\emph{Proof:}
In what follows, without loss of generality, let $n = 1$ for simplicity.
To prove the assertion in theorem \ref{the2}, we implement the following orthogonal transformation from Lemma \ref{lem.wu}:
\begin{subequations}\label{zj}
	\begin{align} \tilde {s }=&{\mathrm{ col}}({\tilde {s }_{1}},{\tilde {s }_{2}}) = {\begin{bmatrix} r& R \end{bmatrix}^{T}} s \\ \tilde {w }=&{\mathrm{ col}}({\tilde {w }_{1}},{\tilde {w }_{2}}) = {\begin{bmatrix} r& R \end{bmatrix}^{T}}w \\ \tilde {s }^{*}=&{\mathrm{ col}}({\tilde {s }_{1}}^{*},{\tilde {s }_{2}}^{*}) = {\begin{bmatrix} r& R \end{bmatrix}^{T}} s ^{*} \\ \tilde {w }^{*}=&{\mathrm{ col}}({\tilde {w }_{1}}^{*},{\tilde {w }_{2}}^{*}) = {\begin{bmatrix} r& R \end{bmatrix}^{T}}w ^{*}\end{align}
\end{subequations}
where $
\tilde{s}_1, \tilde{w}_1, \tilde{s}_{1}^{*}, \tilde{w}_{1}^{*}\in \mathbb{R}$ and $ \tilde{s}_2, \tilde{w}_2, \tilde{s}_{2}^{*}, \tilde{w}_{2}^{*}\in \mathbb{R}^{N-1}.$
As a consequence, one can equivalently rewrite (\ref{al11}) as
\begin{align} \label{zal1}
	\begin{cases}                    %%algo1   ?? 1 
		\dot{x}\in y-x-\partial f\left( y \right) +\left[ r\ R \right] s \\
		\dot{\tilde{s}}_1=-k_1r^T\left( y-d \right) \\
		\dot{\tilde{s}}_2=k_1\left( \tilde{w}_2-R^T\left( y-d \right) \right) -k_2R^TLR\tilde{s}_2\\
		\dot{\tilde{w}}_1=0\\
		\dot{\tilde{w}}_2=k_3R^TL\left( y-d \right) -k_3R^TLR\tilde{w}_2\\
		y=P_{\varOmega}\left( x \right). 
	\end{cases}
\end{align}
To proceed, we only need to analyze the convergence of (\ref{zal1}).

Consider a Lyapunov function candidate as 
\begin{align} \label{v1}  %%%the1 ?? 
	V_1=~&\frac{k_1}{2}\left( \lVert x-P_{\varOmega}\left( x^* \right) \rVert ^2-\lVert x-P_{\varOmega}\left( x \right) \rVert ^2 \right) \nonumber \\  
	&+\frac{1}{2}\lVert \tilde{s }_1-\tilde{s}_{1}^{*} \rVert ^2+\frac{1}{2}\lVert \tilde{s }_2-\tilde{s}_{2}^{*} \rVert ^2 \nonumber \\ 
	&+\frac{1}{2k_3}\lVert \tilde{w }_2-\tilde{w}_{2}^{*} \rVert ^2
\end{align} 
where  $\left( x^*,\tilde{s}_{1}^{*},\tilde{s}_{2}^{*},\tilde{w}_{2}^{*} \right)$ is an equilibrium point of (\ref{m}), and $k_1, k_3$ satisfies \eqref{can1} .
With reference to Lemma \ref{lemproj2}, one can obtain that
\begin{align} \label{zh}
	V_1\ge~& \, \frac{k_1}{2}\lVert y-y^* \rVert ^2 +\frac{1}{2}\lVert \tilde{s} _1-\tilde{s}_{1}^{*} \rVert ^2 \nonumber\\
	&+\frac{1}{2k_3}\lVert \tilde{ w }_2-\tilde{w}_{2}^{*} \rVert ^2+\frac{1}{2}\lVert \tilde{s }_2-\tilde{s}_{2}^{*} \rVert ^2.
\end{align}

The set-valued Lie derivative of $V_1$ along (\ref{zal1}) is expressed as $\mathcal{L}_{\left( \ref{zal1} \right)}V_1$. For any $\zeta _1 \in \mathcal{L}_{\left( \ref{zal1} \right)}V_1$, there exist $\gamma \in \partial f\left( y \right)$ and
$\gamma^{*} \in \partial f\left( y^{*} \right)$  such that 
\begin{align} \label{vdao1}  %%%the1 ?? 
	\zeta _1=~&k_1\lVert y-y^* \rVert ^2-k_1\left( y-y^* \right) ^T\left( x-x^* \right) \nonumber \\ 
	&+k_1\left( \left[ r\,\,R \right] \tilde{s}-\left[ r\,\,R \right] \tilde{s}^* \right) ^T\left( y-y^* \right) \nonumber \\ 
	&-k_1\left( \gamma -\gamma ^* \right) ^T\left( y-y^* \right) \nonumber \\ 
	&-k_1\left( \tilde{s}_2-\tilde{s}_{2}^{*} \right) ^TR^TSym\left( L \right) R\left( \tilde{s}_2-\tilde{s}_{2}^{*} \right) \nonumber \\ 
	&-k_1\left( \tilde{s}_2-\tilde{s}_{2}^{*} \right) ^TR^T\left( y-y^* \right) \nonumber \\ 
	&+k_1\left( \tilde{s}_2-\tilde{s}_{2}^{*} \right) ^T\left( \tilde{w}_2-\tilde{w}_{2}^{*} \right) \nonumber \\ 
	&-\left( \tilde{w}_2-\tilde{w}_{2}^{*} \right) ^TR^TL\left( y-y^* \right) \nonumber \\ 
	&-\left( \tilde{w}_2-\tilde{w}_{2}^{*} \right) ^TR^TSym\left( L \right) R\left( \tilde{w}_2-\tilde{w}_{2}^{*} \right) \nonumber \\ 
	&-k_1\left( \tilde{s}_1-\tilde{s}_{1}^{*} \right) ^Tr^T\left( y-y^* \right) \nonumber \\ 
	=~&k_1\lVert y-y^* \rVert ^2-k_1\left( y-y^* \right) ^T\left( x-x^* \right) \nonumber \\ 
	&-k_1\left( \gamma -\gamma ^* \right) ^T\left( y-y^* \right) \nonumber \\ 
	&-k_2\left( \tilde{s}_2-\tilde{s}_{2}^{*} \right) ^TR^TSym\left( L \right) R\left( \tilde{s}_2-\tilde{s}_{2}^{*} \right) \nonumber \\ 
	&+k_1\left( \tilde{s}_2-\tilde{s}_{2}^{*} \right) ^T\left( \tilde{w}_2-\tilde{w}_{2}^{*} \right) \nonumber \\ 
	&-\left( \tilde{w}_2-\tilde{w}_{2}^{*} \right) ^TR^T L \left( y-y^* \right) \nonumber \\ 
	&-\left( \tilde{w}_2-\tilde{w}_{2}^{*} \right) ^TR^TSym\left( L \right) R\left( \tilde{w}_2-\tilde{w}_{2}^{*} \right)
\end{align} 
where the second equality follows from (\ref{zj}).

With reference to the property of projection given in Lemma \ref{lemproj1} and the strongly convexity of cost function, one can obtain that
\begin{align}   \label{18}
	-\left( y-y^* \right) ^T\left( x-x^* \right) +\lVert y-y^* \rVert ^2\le 0,
\end{align} 
\begin{align}   \label{19}
	-k_1\left( \gamma -\gamma ^* \right) ^T\left( y-y^* \right) \le -k_1\omega \lVert y-y^* \rVert ^2.
\end{align} 
Applying the  Yong's inequality, it yields  
\begin{align}   \label{20}
	k_1\left( \tilde{s}_2-\tilde{s}_{2}^{*} \right) ^T\left( \tilde{w}_2-\tilde{w}_{2}^{*} \right) \le 
	& \frac{k_{1}^{2}}{\hat{\lambda}_2}\lVert s _2-\tilde{s}_{2}^{*} \rVert ^2 \nonumber \\ 
	&+\frac{\hat{\lambda}_2}{4}\lVert w _2-\tilde{w}_{2}^{*} \rVert ^2,
\end{align} 
\begin{align}   \label{21}
	-\left( \tilde{w}_2-\tilde{w}_{2}^{*} \right) ^TR^TL\left( y-y^* \right) 
	\le&  \frac{\hat{\lambda}_2}{4}\lVert w _2-\tilde{w}_{2}^{*} \rVert ^2\nonumber \\  &+\frac{\lVert L \rVert ^2}{\hat{\lambda}_2}\lVert \left( y-y^* \right) \rVert ^2.
\end{align} 
Besides, equipped with the weight-balanced and strongly connected graph, it can be verified that
\begin{align}   \label{22}
	&-k_2\left( \tilde{s}_2-\tilde{s}_{2}^{*} \right) ^TR^TSym\left( L \right) R\left( \tilde{s}_2-\tilde{s}_{2}^{*} \right) \nonumber \\ 
	&\le -k_2\hat{\lambda}_2\lVert s _2-\tilde{s}_{2}^{*} \rVert ^2,
\end{align} 
\begin{align}   \label{23}
	&-\left( \tilde{w}_2-\tilde{w}_{2}^{*} \right) ^TR^TSym\left( L \right) R\left( \tilde{w}_2-\tilde{w}_{2}^{*} \right) \nonumber \\ 
	&\le -\hat{\lambda}_2\lVert w _2-\tilde{w}_{2}^{*} \rVert ^2.
\end{align} 
Consequently, in combination with (\ref{vdao1})-(\ref{23}),  direct calculation yields

\begin{align}   
	\zeta _1\le&-\tau_1
\end{align} 
where 
\begin{align}   
\tau_1=& \left( k_1\omega -\frac{\lVert L \rVert ^2}{\hat{\lambda}_2} \right) \lVert y-y^* \rVert ^2 \nonumber \\ 
	&+\left( k_2\hat{\lambda}_2-\frac{k_{1}^{2}}{\hat{\lambda}_2} \right) \lVert s _2-\tilde{s}_{2}^{*} \rVert ^2 \nonumber \\ 
	&+\frac{\hat{\lambda}_2}{2}\lVert \tilde{w} _2-\tilde{w}_{2}^{*} \rVert ^2\ge 0.  
\end{align} 
It follows from the arbitrariness of $\zeta _1$, that 
\begin{align}\label{ta}   
	\max \mathcal{L}_{(\ref{zal1})}V_1\le -\tau_1\le 0.
\end{align} 

Combining \eqref{ta} and \eqref{zh}, one can arrive at the boundedness of $\left( y\left( t \right) ,\tilde{s}_1\left( t \right) ,\tilde{s}_2\left( t \right) ,\tilde{w}_2\left( t \right) \right) $  for any $t \ge 0$.
Due to the compactness of $\partial f\left(y\right)$, there exists $M>0$ such that 
\begin{align}\label{bou}   
\lVert y\left( t \right) -\gamma +\left[ r\ R \right] s \left( t \right) \rVert \le M,   \forall \gamma  \in \partial f\left( y\left( t \right) \right) , t \ge 0.
\end{align} 
In the sequel, based on the expression of $\dot{x}\left( t \right) $, we show that $x\left( t \right)$  is also bounded.
Define $W\left( x \right) =\frac{1}{2}\lVert x \rVert ^2 $. One can obtain that 
\begin{align}
\mathcal{L}_{(\ref{zal1})}W =\left\{ x^T\left( -x+y-\gamma -\left[ r\ R \right] \tilde{s} \right) :\gamma \in \ \partial f\left( y \right) \right\}.  \nonumber
\end{align}
Subsequently, with $M$ defined in (3), we have
$\max \mathcal{L}_{(\ref{zal1})}W\left(x\left(t \right)\right) \le -\lVert x\left( t \right) \rVert ^2+M\lVert x\left( t \right) \rVert =-2W\left( x\left( t \right) \right) +M\sqrt{2W\left( x\left( t \right) \right)}
$, which can verify the boundness of $x \left(t\right)$.
To proceed, based on Lemma \ref{laser}, the solution of (\ref{al1}) converges to the set $S$ as follows:
\begin{align*}
S=\big\{ \left( x,\tilde{s}_2,\tilde{w}_2 \right) \in ~& \mathbb{R}^N\times \mathbb{R}^{N-1}  \times  \mathbb{R}^{N-1}:\\
& y=y^*,\tilde{s}_2=\tilde{s}_{2}^{*},\tilde{w}_2=\tilde{w}_{2}^{*}  \big\} .
\end{align*}
Thus, it indicates that $\lim _{t\rightarrow \infty}y\left( t \right) =y^*$ and the proof is completed by Theorem \ref{the1}.\hfill $\Box$

\begin{remark}
	Theorem \ref{the2} shows the  effectiveness of algorithm (\ref{al1}) with non-smooth resource allocation under a weight-balanced digraph.  Note that the weight-balanced digraph is a more general assumption than the undirected graph, which may cause the algorithm in \cite{Zeng2016,Xiao06,Lian2021,Lian22021} to not be applicable to the problem solved in Theorem \ref{the2}.
\end{remark}
\begin{theorem}\label{the3}
	For the problem (\ref{op1}) with Assumptions \ref{as1} and \ref{as2}, consider the case where  the local cost functions are strictly convex and the communication topology is a connected and undirected graph.
	 Suppose the initial point $(y(0),x(0),s(0),w(0))$ satisfies $\sum_{i=1}^N{w _i\left( 0 \right) =0_n}
	$, then the algorithm (\ref{al1}) can converge asymptotically to the optimal solution of problem (\ref{op1}).
\end{theorem}
\emph{Proof:}
Similar to the previous proof, we assume $n=1$. Take consider of the following Lyapunov function candidate
\begin{align}            %%%%%The2  ????
	V_2=~&\frac{k_1}{2}\left( \lVert x-P_{\varOmega}\left( x^* \right) \rVert ^2-\lVert x-P_{\varOmega}\left( x \right) \rVert ^2 \right) \nonumber \\ 
	&+\frac{k_3}{2k_2}\lVert \frac{k_1}{k_3}\mathcal{J}_1^{-1}\left( \tilde{w}_2-\tilde{w}_{2}^{*} \right) +\tilde{s}_2-\tilde{s}_{2}^{*} \rVert ^2\nonumber \\
	&+\frac{1}{2}\lVert \tilde{s}_1-\tilde{s}_{1}^{*} \rVert ^2+\frac{1}{2}\lVert \tilde{s} _2-\tilde{s}_{2}^{*} \rVert ^2
\end{align} 
where $\mathcal{J}_1$ is defined in Lemma \ref{lem.wu},  $\left( x^*,\tilde{s}_{1}^{*},\tilde{s}_{2}^{*},\tilde{w}_{2}^{*} \right)$ is an equilibrium point of (\ref{al11}), and $k_1, k_2, k_3 >0$.
With reference to Lemma \ref{lemproj2}, one can obtain that
\begin{align} 
V_2\ge&\frac{k_1}{2}\lVert y-y^* \rVert ^2+\frac{1}{2}\lVert \tilde{s}_1-\tilde{s}_{1}^{*} \rVert ^2+\frac{1}{2}\lVert \tilde{s}_2-\tilde{s}_{2}^{*} \rVert ^2 \nonumber \\ 
&+\frac{k_3}{2k_2}\lVert \frac{k_1}{k_3}\mathcal{J}^{-1}\left( \tilde{w}_2-\tilde{w}_{2}^{*} \right) +\tilde{s}_2-\tilde{s}_{2}^{*} \rVert ^2.
\end{align}

Similarly, for any $\zeta _2 \in \mathcal{L}_{\left( \ref{zal1})\right)}V_2$, there exist $\gamma \in \partial f\left( y \right)$ and $\gamma^{*} \in \partial f\left( y^{*} \right)$  such that
\begin{align}            %%%%%The2  ??????  %%% daoshu2
	\zeta _2=~&k_1\lVert y-y^* \rVert ^2-k_1\left( y-y^* \right) ^T\left( x-x^* \right) \nonumber \\ 
	&+k_1\left( \left[ r\,\,R \right] \tilde{s}-\left[ r\,\,R \right] \tilde{s}^* \right) ^T\left( y-y^* \right) \nonumber \\ 
	&-k_1\left( \gamma -\gamma ^* \right) ^T\left( y-y^* \right) \nonumber \\ 
	&-k_1\left( \tilde{s}_1-\tilde{s}_{1}^{*} \right) ^Tr^T\left( y-y^* \right) \nonumber \\ 
	&-k_1\left( \tilde{s}_2-\tilde{s}_{2}^{*} \right) ^TR^T\left( y-y^* \right) \nonumber \\ 
	&-(k_2+k_3)\left( \tilde{s}_2-\tilde{s}_{2}^{*} \right) ^TR^TLR\left( \tilde{s}_2-\tilde{s}_{2}^{*} \right).
\end{align} 
To proceed, it is straightforward to calculate that
\begin{align}            % %%The2  zuiyou 
	\zeta _2\le -\tau_2
\end{align} 
where 
\begin{align}  
 \tau_2=&\left( k_2+k_3 \right) \left( \tilde{s}_2-\tilde{s}_{2}^{*} \right) ^TR^TLR\left( \tilde{s}_2-\tilde{s}_{2}^{*} \right) \nonumber \\ 
	&+\left( \gamma -\gamma ^* \right) ^T\left( y-y^* \right) \ge 0.
\end{align} 
Then, form the arbitrainess of $\zeta _2$,  it follows that $\max \mathcal{L}_{(\ref{zal1})}V_2\le -\tau_2\le 0.$
%Similarly to the argument in the proof of Theorem \ref{the2},
In a similar way to the arguments shown in the demonstration of Theorem \ref{the2},
 it is verified that $
y\left( t \right) ,\tilde{s}_1\left( t \right) ,\tilde{s}_2\left( t \right) ,\tilde{w}_2\left( t \right) 
$ and $x(t)$ are bounded.  Then in light of Lemma \ref{laser} the solution of (\ref{al11})  is convergent to the largest weakly positively invariant set contained in $Q$, where 
\begin{align*}
	Q=\big\{ \left( x,\tilde{s},\tilde{w}_2 \right) \in ~& \mathbb{R}^N\times \mathbb{R}^N\times \mathbb{R}^{N-1}: \\&\ 0\in \max \mathcal{L}_{\left(\ref{zal1} \right)}V_2\left( x,\tilde{s},\tilde{w}_2 \right) \big\} .
\end{align*}
Furthermore, because of the strict convexity of $f_i$,
one has $
\left( \gamma -\gamma ^* \right) ^T\left( y-y^* \right) >0,\ \forall \ y\ne y^*.
$
Consequently,  $\lim _{t\rightarrow \infty}y\left( t \right) =y^*$.\hfill $\Box$

\begin{remark}
Theorem \ref{the3} shows that algorithm (\ref{al1}) can deal with the problem of non-smooth resource allocation on undirected connected graphs in a fully distributed manner, which means that we do not need to estimate the range of parameters $k_1$ and $k_2$ through additional calculations.
\end{remark}
\begin{remark}
 It is worthwhile to point out that the selecetion of parameters is covering a wide range while it is necessary in \cite{Deng2020} that $k_1=k_2=1$. It is indicated that we can achieve different convergence rates by choosing the appropriate parameters.
\end{remark}

\subsection{Convergence Analysis Of Algorithm \eqref{al2}}\label{C}
In this section, the property corresponding to the equilibrium point of  (\ref{al2}) is first decribed in Theorem \ref{the4}.
Then the non-smooth analysis and the Lyapunov functional theory are employed to demonstrate the convergence of (\ref{al2}) in Theorem \ref{the5}.

Let
$x=col( x_1,x_2,\ldots ,x_N ),\; 
s =col( s _1,s _2,\ldots ,s _N ), \;	
w =col( w _1, w _2,\ldots ,w _N ),\;
d=col( d_1,d_2,\ldots ,d_N ), \;
\varOmega =\varOmega _1 \times \varOmega _2 \times \ldots \times \varOmega _N , \;
\partial f( y ) =col( \partial f_1( y_1 ) ,\partial f_2( y_2 ) , \ldots ,\partial f_N( y_N ) ), \;
y=col( y_1,y_2,\ldots ,y_N ) 
$.

Obviously, algorithm (\ref{al2}) amounts to the following campact form
\begin{align} \label{al22}
	\begin{cases}            %%%%algo2  ????
		\dot{x}\in y-x-\partial f\left( y \right) +s \\
		\dot{s}= k_1\left( \left( L\otimes I_n \right) w -y+d \right) -k_2\left( L\otimes I_n \right) s \\
		\dot{w}=-k_3\left( L\otimes I_n \right) \left( \left( L\otimes I_n \right) w -y+d \right) \\
		y=P_{\varOmega}\left( x \right) .
	\end{cases}
\end{align}

By virtue of Lemma \ref{map} and Assumption \ref{as2}, it is observed that the system (\ref{al22}) exists a solution.

The result given in the following is concerning the equilibrium point of (\ref{al22}).

\begin{theorem}\label{the4}
	For the nonsmooth resource allocation problem (\ref{op1}), consider the case where the communication topology is a connected and undirected graph.
 If Assumptions \ref{as1} and \ref{as2} hold, then $(y^*,x^*,s^*,w^*)$ is an equilibrium point of \eqref{al22}, if and only if $y^*$ is an optimal solution of the problem (\ref{op1}).
\end{theorem}
 \emph{Proof:}  
 1) Let $(y^*,x^*,s^*,w^*)$ be an  equilibrium point of (\ref{al22}), then one has \begin{subequations}
	\begin{align}\label{Th2a}    %%%%  algo2???
		0_{Nn}\in ~&y^{\ast}-x^{\ast}-\partial f\left( y^{\ast} \right) +s ^{\ast}\\ \label{Th2b}
		0_{Nn}=~&k_1\left( \left( L\otimes I_n \right) w ^{\ast}-y^{\ast}+d \right) -k_2\left( L\otimes I_n \right) s ^{\ast}\\ \label{Th2c}
		0_{Nn}=~&-k_3\left( L\otimes I_n \right) \left( \left( L\otimes I_n \right) w ^{\ast}-y^{\ast}+d \right) \\ \label{Th2d}
		y^{\ast}=~&P_{\varOmega}\left( x^{\ast} \right) .
	\end{align}
\end{subequations}

Firstly, according to  (\ref{Th2c}), there exists $\theta_2 \in \mathbb{R}^n$ such that $
\left( L\otimes I_n \right) w ^{\ast}-y^{\ast}+d=1_N\otimes \theta _2 $ due to the connectedness of the undirected graphs. 
Thus, (\ref{Th2b}) indicates that $k_2\left( L\otimes I_n \right) s ^*=k_1\left( 1_N\otimes \theta _2 \right)$ and $k_2\left( L^TL\otimes I_n \right) s ^*=k_1\left( L^T1_N\otimes \theta _2 \right) =0_{Nn}$.
Based on Lemma \ref{lem.dac}, it follows that $\left( L\otimes I_n \right) s ^{\ast}=0_{Nn}$, i.e. $s _{i}^{*}=s _{j}^{*}, \forall \, i,j \in \{1,\dots,N\}$.

Next, form (\ref{Th2b}), it can be calculated that $k_1\left( 1_{N}^{T}\otimes I_n \right) \left( -y^*+d \right) =k_2\left( 1_{N}^{T}L\otimes I_n \right) s ^*-k_1\left( 1_{N}^{T}L\otimes I_n \right) w ^*=0_{n}$,  which signifies that $\sum_{i=1}^N{y_{i}^{*}=}\sum_{i=1}^N{d_i}$.
%k_1\left( 1_{N}^{T}\otimes I_n \right) \left( \left( L\otimes I_n \right) w ^*-y^*+d \right) =

Finally, (\ref{Th2a}) and (\ref{Th2d}) imply that $y^*=P_{\varOmega}\left( y^*-\partial f\left( y^* \right) +s ^* \right)$ which is identical to $0_n\in \partial f_i\left( y_{i}^{*} \right) -s _{i}^{*}+N_{\varOmega _i}\left( y_{i}^{*} \right)$, $\forall \in \{1,\ldots,N\}$.

Combined with the above discussions, it yields that $y^*$ is an optimal solution of the problem \eqref{op1} by Lemma  \ref{kkt}.

  2) If $y^*$ is an optimal solution of the problem (\ref{op1}), then there exists $s _{i}^{*}\in \mathbb{R}^n$  satisfying (\ref{zui}).

Thus, we can take $s ^*=col\left\{ s _{1}^{*},\dots,s _{N}^{*} \right\}$ and  claim that $\left( L\otimes I_n \right) s ^*=0_{Nn}$ by (\ref{zui}). Based on (\ref{zui}), setting $x^*=col\left\{ x_{1}^{*},\dots,x_{N}^{*} \right\}$ where $x_{i}^{*}\in y_{i}^{*}-\partial f_i\left( y_{i}^{*} \right) +s _{i}^{*}$, it is easy to verify that (\ref{Th2a}) and (\ref{Th2b}) hold.

In the sequel, we illustrate there exists $w ^*\in \mathbb{R}^n$ such that $\left( L\otimes I_n \right) w ^*=y^*-d$, which can give rise to the satisfaction of (\ref{Th2b}) and (\ref{Th2c}) without much effort. It follows from the connectedness of the undirected graphs that $1_N\otimes \xi  \in ker\left( L\otimes I_n \right)$ for any $\xi \in \mathbb{R}^n$. Meanwhile, one can get $\left( y^*-d \right) ^T\left( 1_N\otimes \xi  \right) =0$ due to (\ref{zui}). As a consequence, $y^*-d\in range\left( L\otimes I_n \right)$ by noting that $\mathbb{R}^{Nn}$ can be orthogonally decomposed by $ker\left( L\otimes I_n \right)$ and $range\left( L\otimes I_n \right)$ (see \cite{s1993The}).
Hence, there exists $w ^*\in \mathbb{R}^{Nn}$ such that $(L\otimes I_n )w^*=y^*-d$.

With reference to the above analysis, it can be concluded that $\left( y^*,x^*,s ^*,w ^* \right)$ is an equilibrium point of (\ref{al22}). \hfill $\Box$

Next, the result about the convergence of (\ref{al2}) over an undirected graph is presented.
	\begin{theorem}\label{the5}
		For the problem (\ref{op1}) with Assumptions \ref{as1} and \ref{as2}, consider the case where  the local cost functions are $\omega$-strongly convex and the communication topology is a connected undirected graph.
		 The algorithm (\ref{al2}) can converge asymptotically to the optimal solution of problem (\ref{op1}).
\end{theorem}
\emph{Proof:}
In the sequel, by assigning $n=1$, we examine the following equivalent formulation of algorithm \eqref{al22} obtained by the orthogonal transformation \eqref{zj}:
\begin{align} \label{zal2}
	\begin{cases}               %%algo2  ?? 2
		\dot{x}\in y-x-\partial f\left( y \right) +\left[ r\ R \right] s \\
		\dot{\tilde{s}}_1=-r^T\left( y-d \right) \\
		\dot{\tilde{s}}_2=k_1\left( R^TLR\tilde{w}_2-R^T\left( y-d \right) \right) -k_2R^TLR\tilde{s}_2\\
		\dot{\tilde{w}}_1=0\\
		\dot{\tilde{w}}_2=k_3R^TL\left( y-d \right) -k_3R^TL^2R\tilde{w}_2\\
		y=P_{\varOmega}\left( x \right). 
	\end{cases}
\end{align}
Therefore, we only need to discuss the convergence of (\ref{zal2}).

Select the Lyapunov function candidate as
\begin{align}   
	V_3=&\frac{k_1}{2}\left( \lVert x-P_{\varOmega}\left( x^* \right) \rVert ^2-\lVert x-P_{\varOmega}\left( x \right) \rVert ^2 \right) \nonumber \\ 
	&+\frac{1}{2}\lVert s _1-\tilde{s}_{1}^{*} \rVert ^2+\frac{1}{2}\lVert s _2-\tilde{s}_{2}^{*} \rVert ^2 \nonumber \\ 
	&+\frac{1}{2k_3}\lVert w _2-\tilde{w}_{2}^{*} \rVert ^2
\end{align} 
where  $\left( x^*,\tilde{s}_{1}^{*},\tilde{s}_{2}^{*},\tilde{w}_{2}^{*} \right)$ is an equilibrium point of (\ref{al22}), and $k_1, k_3 >0$.
With reference to Lemma \ref{lemproj2}, one can obtain that
\begin{align} 
V_3\ge&\frac{k_1}{2}\lVert y-y^* \rVert ^2+\frac{1}{2}\lVert \tilde{s}_1-\tilde{s}_{1}^{*} \rVert ^2\nonumber \\ 
&+\frac{1}{2}\lVert \tilde{s}_2-\tilde{s}_{2}^{*} \rVert ^2
+\frac{1}{2k_3}\lVert \tilde{w}_2-\tilde{w}_{2}^{*} \rVert ^2.
\end{align} 

Similar to the previous arguments, consider the ser-valued Lie derivative of $V_3$ with respect to (\ref{zal2}). For any $\zeta _1 \in \mathcal{L}_{\left( \ref{zal2} \right)}V_1$, there exist $\gamma \in \partial f\left( y \right)$ and
$\gamma^{*} \in \partial f\left( y^{*} \right)$  such that 
\begin{align}   
	\zeta _3=~&k_1\lVert y-y^* \rVert ^2-k_2\left( y-y^* \right) ^T\left( x-x^* \right)   \nonumber \\ 
	&+k_1\left( \left[ r\,\,R \right] \tilde{s}-\left[ r\,\,R \right] \tilde{s}^* \right) ^T\left( y-y^* \right)  \nonumber \\ 
	&-k_1\left( \gamma -\gamma ^* \right) ^T\left( y-y^* \right)  \nonumber \\ 
	&-k_2\left( \tilde{s}_2-\tilde{s}_{2}^{*} \right) ^TR^T L  R\left( \tilde{s}_2-\tilde{s}_{2}^{*} \right)  \nonumber \\ 
	&-k_1\left( \tilde{s}_2-\tilde{s}_{2}^{*} \right) ^TR^T\left( y-y^* \right)  \nonumber \\ 
	&+k_1\left( \tilde{s}_2-\tilde{s}_{2}^{*} \right) ^TR^T L  R\left( \tilde{w}_2-\tilde{w}_{2}^{*} \right)  \nonumber \\ 
	&-\left( \tilde{w}_2-\tilde{w}_{2}^{*} \right) ^TR^TL\left( y-y^* \right)  \nonumber \\ 
	&-\left( \tilde{w}_2-\tilde{w}_{2}^{*} \right) ^TR^T L^2  R\left( \tilde{w}_2-\tilde{w}_{2}^{*} \right)  \nonumber \\ 
	&-k_1\left( \tilde{s}_1-\tilde{s}_{1}^{*} \right) ^Tr^T\left( y-y^* \right)  \nonumber \\ 
	=~&k_1\lVert y-y^* \rVert ^2-k_2\left( y-y^* \right) ^T\left( x-x^* \right)  \nonumber \\ 
	&-k_1\left( \gamma -\gamma ^* \right) ^T\left( y-y^* \right)  \nonumber \\ 
	&-k_2\left( \tilde{s}_2-\tilde{s}_{2}^{*} \right) ^TR^T L  R\left( \tilde{s}_2-\tilde{s}_{2}^{*} \right)  \nonumber \\ 
	&+k_1\left( \tilde{s}_2-\tilde{s}_{2}^{*} \right) ^TR^T L R\left( \tilde{w}_2-\tilde{w}_{2}^{*} \right)  \nonumber \\ 
	&-\left( \tilde{w}_2-\tilde{w}_{2}^{*} \right) ^TR^TL \left( y-y^* \right)  \nonumber \\ 
	&-\left( \tilde{w}_2-\tilde{w}_{2}^{*} \right) ^TR^T L^2 R\left( \tilde{w}_2-\tilde{w}_{2}^{*} \right). 
\end{align} 
According to Yong's inequality, it implies   
\begin{subequations}\label{39}
	\begin{align}    %%% %%%???????
	&k_1\left( \tilde{s}_2-\tilde{s}_{2}^{*} \right) ^TR^TL  R\left( \tilde{w}_2-\tilde{w}_{2}^{*} \right)  \nonumber \\ 
& \le \frac{k_{1}^{2}\lVert L \rVert ^2}{{\lambda}_2^2}\lVert \tilde{s}_2-\tilde{s}_{2}^{*} \rVert ^2+\frac{{\lambda}_2^2}{4}\lVert \tilde{w}_2-\tilde{w}_{2}^{*} \rVert ^2,\\
&-\left( \tilde{w}_2-\tilde{w}_{2}^{*} \right) ^TR^TL\left( y-y^* \right)   \nonumber \\
&\le \frac{{\lambda}_2^2}{4}\lVert \tilde{w}_2-\tilde{w}_{2}^{*} \rVert ^2+\frac{\lVert L \rVert ^2}{{\lambda}_2^2}\lVert \left( y-y^* \right) \rVert ^2.
	\end{align} 
\end{subequations}
Utilized the fact that the connected graph is undirected and connected, one can obtain that
\begin{subequations}\label{40}  
\begin{align}     %%% %%%???????
-k_2\left( \tilde{s}_2-\tilde{s}_{2}^{*} \right) ^TR^T L R\left( \tilde{s}_2-\tilde{s}_{2}^{*} \right) &\le -k_2{\lambda}_2\lVert \tilde{s}_2-\tilde{s}_{2}^{*} \rVert ^2, \\
-\left( \tilde{w}_2-\tilde{w}_{2}^{*} \right) ^TR^TL^2  R\left( \tilde{w}_2-\tilde{w}_{2}^{*} \right)& \le -{\lambda}_2^2\lVert \tilde{w}_2-\tilde{w}_{2}^{*} \rVert ^2.
\end{align} 
\end{subequations}

As a result, combined with (\ref{18}) (\ref{19}) (\ref{39}) (\ref{40}), $\zeta _3$ can be ecaluated as    
\begin{align}\label{tao}       %%????
	\zeta _3\le &-\left( k_1\omega -\frac{\lVert L \rVert ^2}{{\lambda}_2^2} \right) \lVert y-y^* \rVert ^2\nonumber \\ 
	&-\left( k_2{\lambda}_2-\frac{k_{1}^{2}\lVert L \rVert ^2}{{\lambda}_2^2} \right) \lVert \tilde{s} _2-\tilde{s}_{2}^{*} \rVert ^2\nonumber \\ 
	&-\frac{{\lambda}_2^2}{2}\lVert \tilde{w} _2-\tilde{w}_{2}^{*} \rVert ^2.
\end{align} 

In view of the above inequalities (\ref{tao}), the remainder of this proof can be derived by similar argument in the proof of Theorem \ref{the2} and, thus, can be omitted.\hfill $\Box$

\begin{remark}
	Summarizing the above analysis, it is obvious that the auxiliary variable $\mu_i$ in algorithm (\ref{al1}) needs to satisfy certain initial conditions, while algorithm (\ref{al2}) can implement by an initialization-free way.  In addition, compared with algorithm (\ref{al1}), the dynamics of algorithm (\ref{al2}) are more complicated and have higher communication costs.  However, some practical application scenarios are not easy to fulfill the initial conditions in algorithm (\ref{al1}).  In this instance, it is worthwhile to sacrifice a certain communication cost to realize the algorithm without initialization.
\end{remark}

\section{Numerical Example}\label{section5}
In this section, we provide two numerical examples to illustrate the utility and preformance of the proposed algorithms.

\textbf{Example 1:}  In this subsection, inspired by \cite{Deng2018}, we verify the effectiveness of the algorithms through an economic dispatch problem in a smart grid. Specifically, we consider a grid composed of four generators with the interaction network described by $\mathcal{G}_1$ or $\mathcal{G}_2$ shown in Fig.\ref{fig.g}. The cost function of each generator in $M\$$ takes the form of a nonsmooth quadratic function as follows:
\begin{figure}[!t]
	\centering 	
	\subfloat[$\mathcal{G}_1$]{{\includegraphics[width=1.2in]{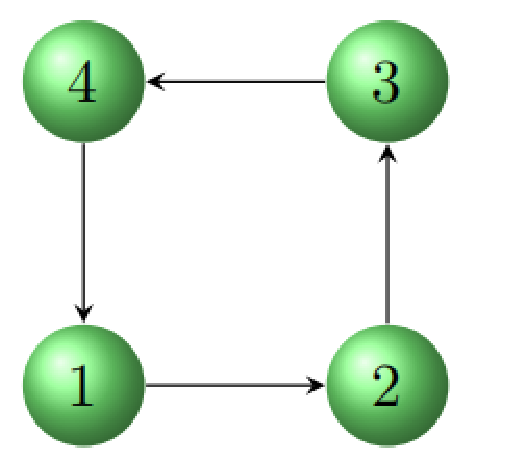} }}	
	\subfloat[$\mathcal{G}_2$]{{\includegraphics[width=1.2in]{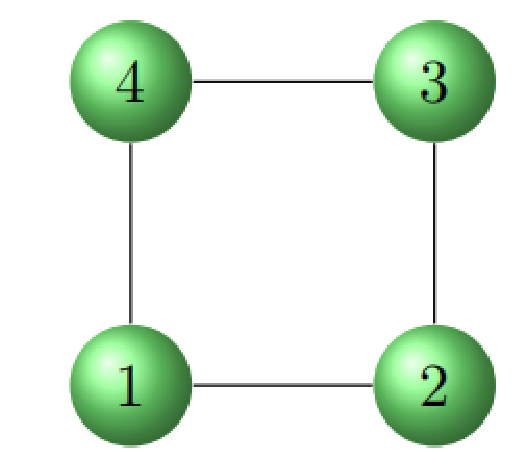} }}%
	\caption{Communication graph of the network: (a) directed and weight-balanced; (b) undirected.}
	\label{fig.g}%
\end{figure}
\begin{align*}
	f_i(p_{Gi})= \alpha_i +\beta_i \|p_{Gi}-35\|+\gamma_i p_{Gi}^2
\end{align*}
where $p_{Gi}$ represents the output power in $MW$ and $\alpha_i, \beta_i $ and $\gamma_i$ denote the system parameters of generator $i$.
Further, for each generator $i$, $p_{d_i}$ is employed to denote the local load demand.

For security, economic and other factors, in practical applications, the output power  is usually specified to have a certain upper boundary $p_{Gi}^{\max}$ and a lower boundary $p_{Gi}^{\min}$, that is  $p_{Gi}^{\min} \le p_{Gi} \le p_{Gi}^{\max}$.
Specifically, for $i$th generator, the parameters in cost functions, the upper and lower bounds of output power, and the local load demand are listed in Table \ref{table:parameters}.
\begin{table}
	\begin{center}
		\caption{Parameters setting \cite{Deng2018}}\label{table:parameters}
		\begin{tabular}{ccccccc}
			\hline
			Generator &	$\alpha_i$  &$\beta_i$ & $\gamma_i$& $p_{G_i}^{\min}$ & $p_{G_i}^{\max}$ & $p_{d_i}$
			\\ \hline
			$1$ & $0.5$ & $3$    & $2$ & $20$ & $40$& $45$ \\
			$2$ & $1.5$ & $4$    & $1$ & $25$ & $35$& $40$ \\
			$3$ & $3$ & $5$    & $0.5$ & $35$ & $50$& $25$  \\
			$4$ & $1$ & $2$   & $1.5$ & $25$ & $45$ & $35$ \\
			\hline
		\end{tabular}
	\end{center}
\end{table}

\begin{figure}
	\centering
	\includegraphics[width=9cm]{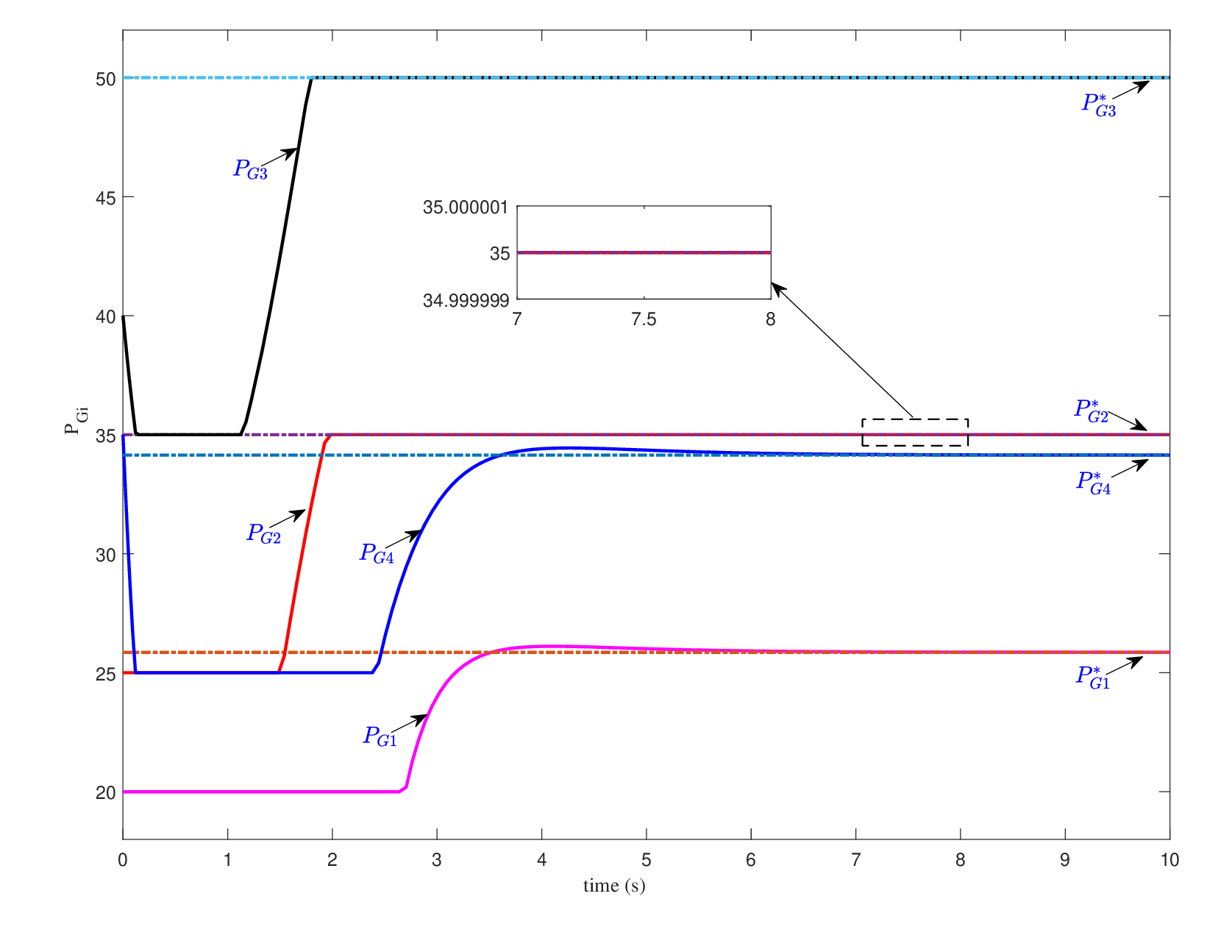}\\
	\caption{Simulation results of Example 1 by  algorithm (\ref{al1}).}\label{fig.1}
\end{figure}

\begin{figure}
	\centering
	\includegraphics[width=9cm]{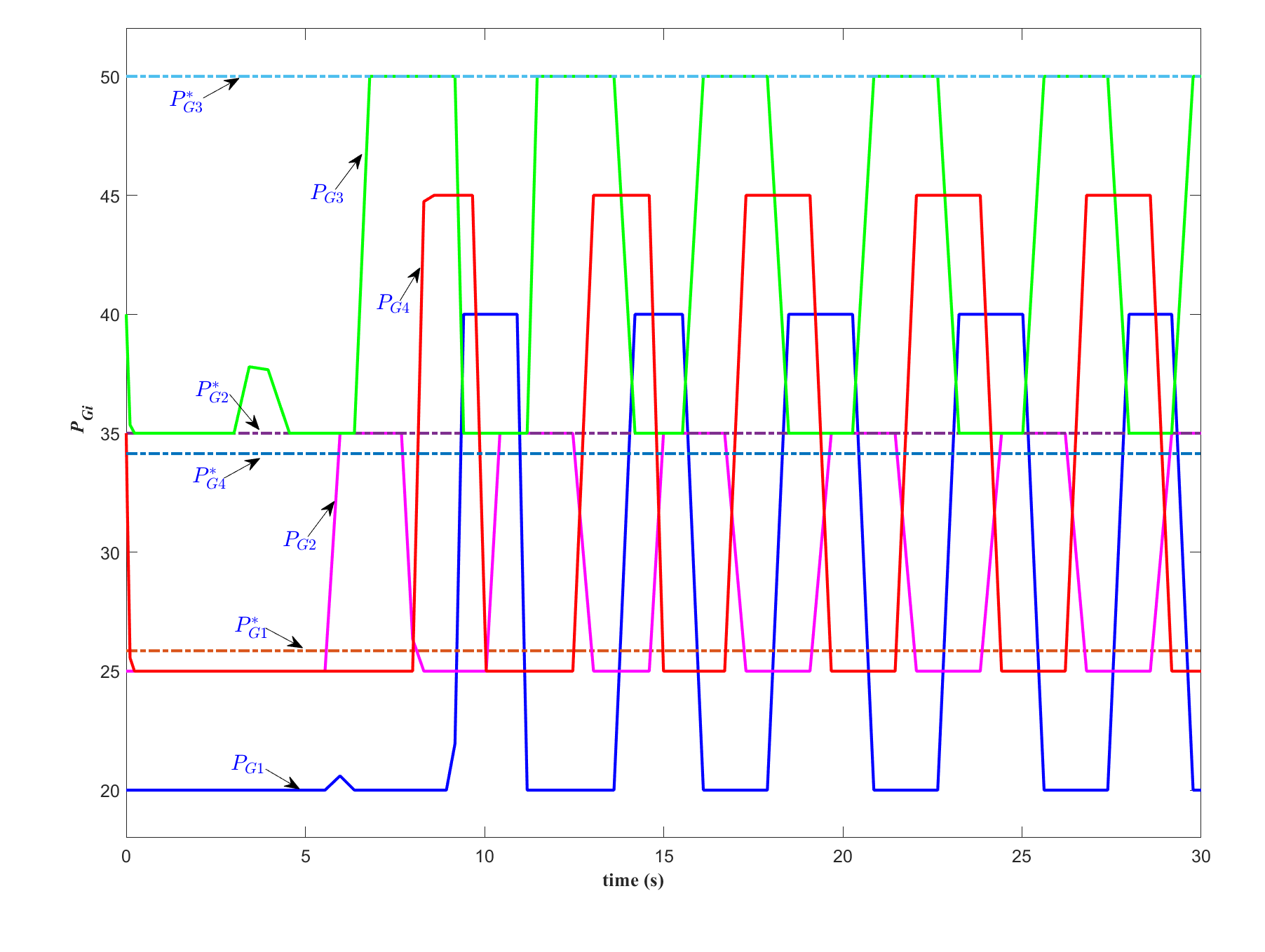}\\
	\caption{Simulation results of Example 1 by the  algorithm proposed in \cite{Zeng2016}.}\label{fig.2}
\end{figure}

By executing algorithm \eqref{al1} under the weight-balanced digraph $\mathcal{G}_1$ shown in Fig.\ref{fig.g}(a), setting the parameters as $k_1=5, k_2=26,$  and $k_3=5$, and assigning auxiliary variables as  $w _i\left( 0 \right) =0, i \in \{1,\ldots,4\} $, we can obtain the simulation results as depicted in Fig.\ref{fig.2}. We use dotted lines to indicate the optimal output power associated with the generator in fig.\ref{fig.g}. It can be observed from Fig.\ref{fig.1} that the optimal value is $
p_G^*=(25.8569, \,35.0000, \,50.0000, \,34.1431)^T$ and the output power $p_i$ of each generator converges to the exact optimal value. Moreover, with a simple calculation, it can be checked that the final output powers satisfy the network resource constraint.

\begin{figure}
	\centering
	\includegraphics[width=9cm]{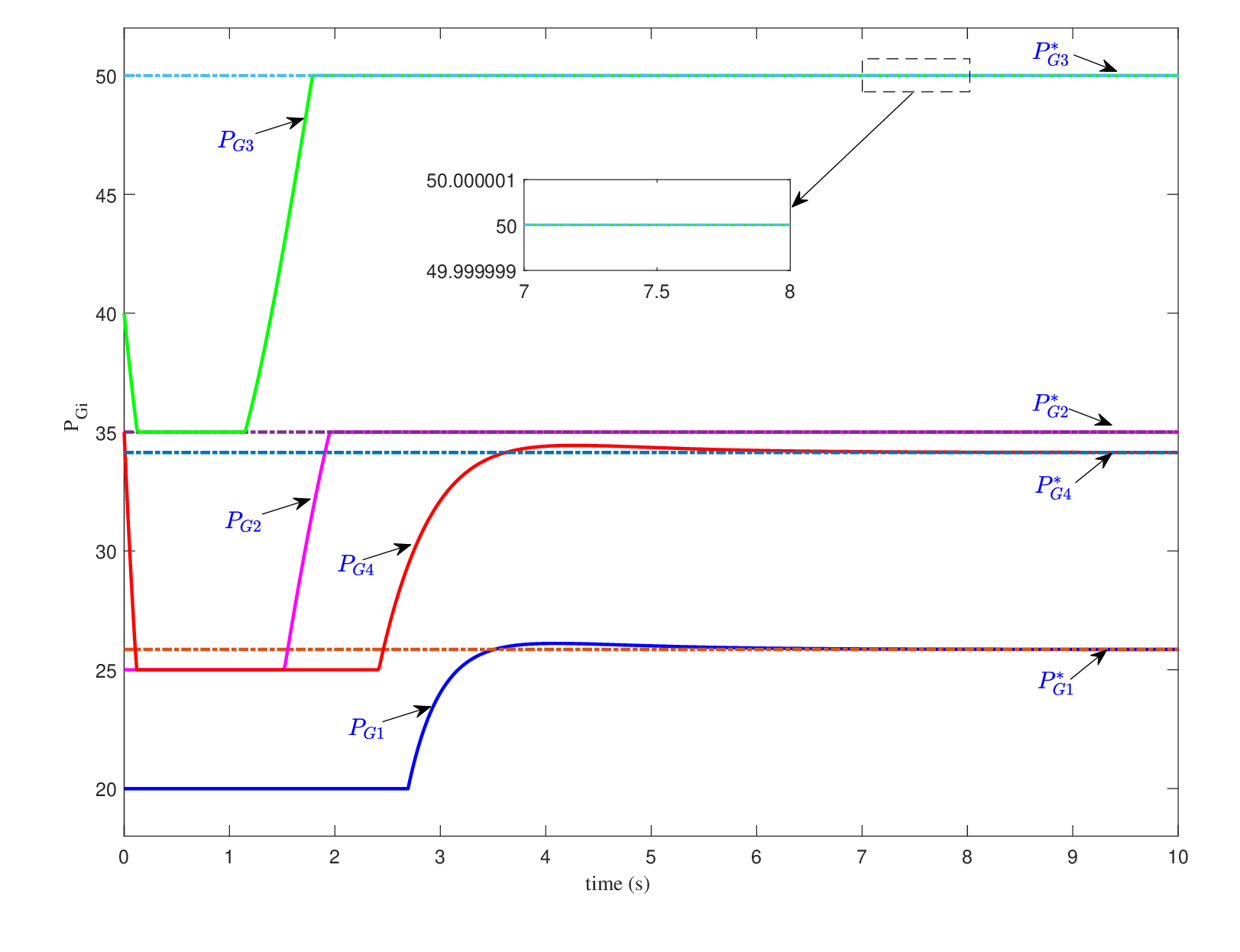}\\
	\caption{Simulation results of Example 1 by  algorithm (\ref{al2}).}\label{fig.3}
\end{figure}

It is worth mentioning that, as illustrated in Fig.\ref{fig.2}, the algorithm in \cite{Zeng2016} cannot be suitable in the case that the communication network is the directed graph $\mathcal{G}_1$ shown in Fig.\ref{fig.g}(a).  As a comparison, algorithm \eqref{al1} developed in this paper overcomes the obstacles caused by a directed graph which means that it can be applied in more general practical scenarios.

Running algorithm \eqref{al2} on the undirected graph $\mathcal{G}_2$  shown in Fig.\ref{fig.g}(b), setting the parameters to $k_1=5,k_2=55,k_3=5$, and configuring the auxiliary variables  as  $w _i\left( 0 \right) =10, i \in \{1,2,3\} $ and $w_4(0)=0$, we can obtain the simulation results presented in Fig.\ref{fig.3}, where the optimal output power is indicated using dotted lines. Compared to Fig.\ref{fig.1}, it is clear to observe that the output power of each generator also converges to the optimal value, even though we did not set the initial values of the auxiliary variables to satisfy specific requirements, thanks to the initialization-free nature of algorithm (\ref{al2}).

\textbf{Example 2:}  Note that cost functions much more complicated than the quadratic function in Example 1 frequently appear in practical engineering. To further exemplify the generalizability of the algorithm developed in this paper, we next consider some more complex cost functions.

We consider the following problem of distributed resource allocation for four agents communicating via the two graphs shown in Fig.\ref{fig.g}, respectively.

The local cost function $f_i (y_i)$, local feasible set constraint $\varOmega_i$, and local load demand $d_i$ for each  agent $i\in \{1,\dots,4\},$ are defined as follows:
\begin{align*}
	f_1\left( y_1 \right) =~&\lVert y_1 \rVert ^2+\lVert y_1-\left[ 2\ 2 \right] ^T \rVert \\
	f_2\left( y_2 \right) =~&\lVert y_2 \rVert ^2+\frac{y_{21}^{2}}{20{y_{21}^{2}+1}}+\frac{y_{22}^{2}}{20{y_{22}^{2}+1}}\\
	f_3\left( y_3 \right)=~&\lVert y_3-\left[ 2\ 3 \right] ^T \rVert ^2\\
	f_4\left( y_4 \right) =~&\ln \left( e^{-0.05y_{41}}+e^{0.05y_{41}} \right) \\
	&+\ln \left( e^{-0.05y_{42}}+e^{0.05y_{42}} \right) +\lVert y_4 \rVert ^2\\
	\varOmega_1=~&\left\{ y_1\in R^2\left|~ \lVert y_1-\left[ 2\ 2 \right] ^T \rVert \le 2 \right. \right\} \\
	\varOmega_2=~&\left\{ y_2\in R^2\left| ~1\le y_{21}\le 2,0\le y_{22}\le 1 \right. \right\} \\
	\varOmega_3=~&\left\{ y_3\in R^2\left| ~y_{31}\ge 0.5,y_{32}\ge 1,y_{31}+y_{32}\le 6 \right. \right\} \\
	\varOmega_4=~&\left\{ y_4\in R^2\left| ~\lVert y_4-\left[ 3\ 5 \right] ^T \rVert \le 2 \right. \right\} 
\end{align*}
where $y_i =(y_{i1}, y_{i2})^T \in \mathbb{R}^2$, 
and  $d_1=(2, \; 1)^T$, $d_2=(2, \; 3)^T$, $d_3=(2, \; 4)^T$, and $d_3=(1, \; 5)^T$, respectively.

\begin{figure}
	\centering
	\includegraphics[width=9cm]{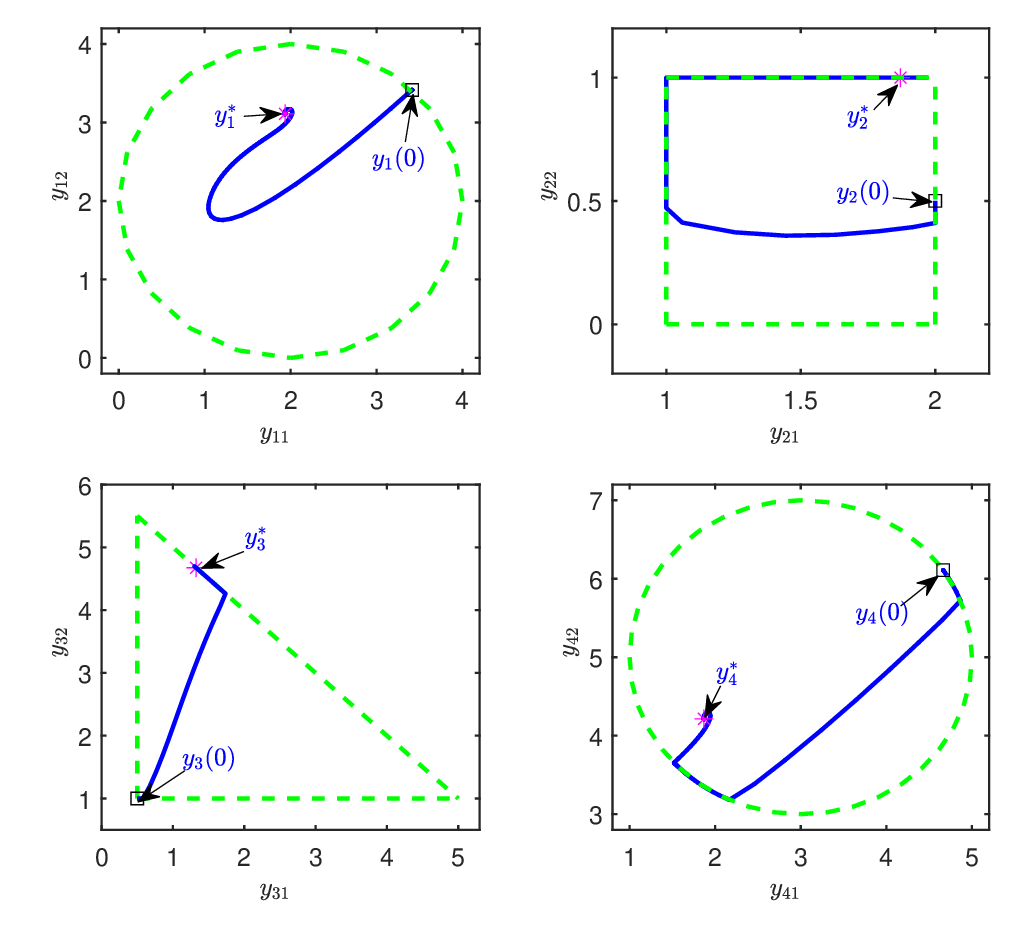}\\
	\caption{Evolutionary trajectories of the decisions of four agents by algorithm (\ref{al1}) under Fig.\ref{fig.g}(a).}\label{fig.4}
\end{figure}

\begin{figure}
	\centering
	\includegraphics[width=9cm]{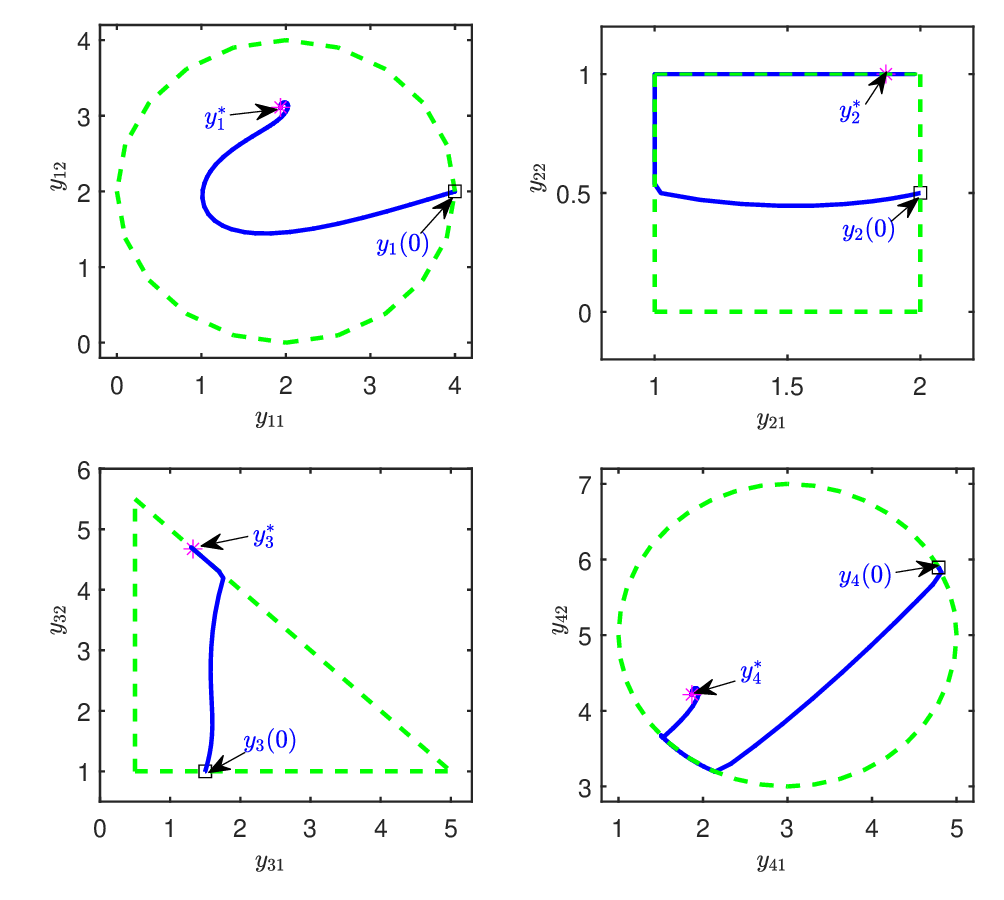}\\
	\caption{Evolutionary trajectories of the decisions of four agents by algorithm (\ref{al1}) under Fig.\ref{fig.g}(b).}\label{fig.5}
\end{figure}

Fig.\ref{fig.4} and Fig.\ref{fig.5} depict the simulation results of algorithm (\ref{al1}) under the weight-balanced digraph shown in Fig.\ref{fig.g}(a) and the undirected graph shown in Fig.\ref{fig.g}(b), respectively, where the green lines signify the local feasible set constraints and the blue lines signify the evolution of the decisions. 
In Fig.\ref{fig.3}, setting the parameters of algorithm (\ref{al1}) as $k_1=5,k_2=26$ and $k_3=5$, it can be noticed that the evolutionary trajectory of the decision for each agent always stays within the local feasible set and finally converges to the optimal solution.

In Fig.\ref{fig.4}, setting the parameters of algorithm \ref{al1} as $k_1=5,k_2=5$ and $k_3=5$, we can observe that the decisions  also effectively converge to the optimal solution. In particular, in Fig.\ref{fig.5}, we set $k_2$ that does not satisfy the condition in (\ref{can1}), yet algorithm (\ref{al1}) is still valid, which shows that we can run algorithm (\ref{al1}) in a fully distributed manner under undirected graphs without additional restrictions on the parameter range, as stated in Theorem \ref{the3}.

\section{Conclusion}\label{section6}

In this paper, the nonsmooth resource allocation problem with heterogeneous constraints depicted by general convex sets is investigated. We develop a novel distributed algorithm via differential inclusion and projected output feedback. It is proved that the algorithm can solve the nonsmooth resource allocation problem on weight-balanced digraphs with strongly convex cost functions. Furthermore, the algorithm is also proved to resolve the problem on undirected graphs in a fully distributed manner with strictly convex cost functions. In addition, a new algorithm is developed to improve the drawback requirement of the initialization of auxiliary variables in the first algorithm. The initialization-free algorithm is proved to address the nonsmooth resource allocation problem on undirected graphs with strongly convex cost functions by the Lyapunov functional theory and the nonsmooth analysis theory. Besides, some simulations are carried out for the effectiveness of proposed algorithms. Future work may focus on the nonsmooth resource allocation problem in the presence of communication delays or external disturbances.

\ifCLASSOPTIONcaptionsoff
  \newpage
\fi

\bibliographystyle{IEEEtran}
\bibliography{resource}

\end{document}